\documentclass[11pt,reqno,a4paper]{amsart}

\usepackage{a4,amssymb}
\usepackage[ansinew]{inputenc}
\allowdisplaybreaks   

\usepackage{hyperref}

\theoremstyle{plain}
\newtheorem{de}{Definition}[section]
\newtheorem{lem}[de]{Lemma}
\newtheorem{prop}[de]{Proposition}

\newtheorem{thm}[de]{Theorem}

\theoremstyle{definition}
\newtheorem{rem}[de]{Remark}

\newcommand{\R}{\mathbb{R}}

\newcommand{\C}{\mathbb{C}}

\newcommand{\D}{\,\textrm{d}}

\newcommand{\trm}{\textrm}

\newcommand{\bea}{\begin{eqnarray*}}
\newcommand{\eea}{\end{eqnarray*}}

\newcommand{\beq}{\begin{equation}}
\newcommand{\eeq}{\end{equation}}

\newcommand{\bprf}{\begin{proof}[\emph{\textbf{Proof.}}]}
\newcommand{\eprf}{\end{proof}}

\newcommand{\ra}{\rightarrow}

\newcommand{\ol}{\overline}

\newcommand{\wt}{\widetilde}
\newcommand{\bs}{\backslash}
\newcommand{\ii}{\textrm{i}}

\newcommand{\im}{\textrm{im}\,}

\newcommand{\OO}{\mathcal{O}}

\newcommand{\uu}{\texttt{u}}
\newcommand{\ue}{\emph{\texttt{u}}}

\newcommand{\rea}{\textrm{Re}\,}
\newcommand{\ima}{\textrm{Im}\,}
\newcommand{\id}{\textrm{id}\,}

\newcommand{\eps}{\varepsilon}

\def\typeout#1{\message{^^J}\message{#1}\message{^^J}}
\newif\ifSRCOK \SRCOKtrue
\newcount\PAGETOP
\newcount\LASTLINE
\global\PAGETOP=1
\global\LASTLINE=-1
\def\EJECT{\SRC\eject}
\def\WinEdt#1{\typeout{:#1}}% WinEdt LOG MODE and INPUT
\gdef\MainFile{\jobname.tex}% ".tex" needed for MiKTeX
\gdef\CurrentInput{\MainFile}
\newcount\INPSP
\global\INPSP=0
\def\SRC{\ifSRCOK%
  \ifnum\inputlineno>\LASTLINE%
    \ifnum\LASTLINE<0%
      \global\PAGETOP=\inputlineno%
    \fi%
    \global\LASTLINE=\inputlineno%
    \ifnum\INPSP=0%
      \ifnum\inputlineno>\PAGETOP%
        
      \fi%
    \else%
      
    \fi%
  \fi%
\fi}
\def\PUSH#1{%
\SRC%
\ifnum\INPSP=0 \global\let\INPSTACKA=\CurrentInput \else%
\ifnum\INPSP=1 \global\let\INPSTACKB=\CurrentInput \else%
\ifnum\INPSP=2 \global\let\INPSTACKC=\CurrentInput \else%
\ifnum\INPSP=3 \global\let\INPSTACKD=\CurrentInput \else%
\ifnum\INPSP=4 \global\let\INPSTACKE=\CurrentInput \else%
\ifnum\INPSP=5 \global\let\INPSTACKF=\CurrentInput \else%
               \global\let\INPSTACKX=\CurrentInput \fi\fi\fi\fi\fi\fi%
\gdef\CurrentInput{#1}%
\WinEdt{<+ \CurrentInput}%
\global\LASTLINE=0%
\ifSRCOK\fi%
\global\advance\INPSP by 1}
\def\POP{%
\ifnum\INPSP>0 \global\advance\INPSP by -1  \fi%
\ifnum\INPSP=0 \global\let\CurrentInput=\INPSTACKA \else%
\ifnum\INPSP=1 \global\let\CurrentInput=\INPSTACKB \else%
\ifnum\INPSP=2 \global\let\CurrentInput=\INPSTACKC \else%
\ifnum\INPSP=3 \global\let\CurrentInput=\INPSTACKD \else%
\ifnum\INPSP=4 \global\let\CurrentInput=\INPSTACKE \else%
\ifnum\INPSP=5 \global\let\CurrentInput=\INPSTACKF \else%
               \global\let\CurrentInput=\INPSTACKX \fi\fi\fi\fi\fi\fi%
\WinEdt{<-}%
\global\LASTLINE=\inputlineno%
\global\advance\LASTLINE by -1%
\SRC}
\def\INPUT#1{\relax}
\let\originalxxxeverypar\everypar
\newtoks\everypar
\originalxxxeverypar{\the\everypar\expandafter\SRC}
\everymath\expandafter{\the\everymath\expandafter\SRC}
\output\expandafter{\expandafter\SRCOKfalse\the\output}
\newif\ifSRCOK \SRCOKtrue
\newcount\PAGETOP
\newcount\LASTLINE
\global\PAGETOP=1
\global\LASTLINE=-1
\gdef\MainFile{\jobname.tex}% ".tex" needed for MiKTeX
\gdef\CurrentInput{\MainFile}
\newcount\INPSP
\global\INPSP=0
\def\EJECT{\SRC\eject}
\def\WinEdt#1{\typeout{:#1}}% WinEdt LOG MODE and INPUT
\def\SRC{\ifSRCOK%
  \ifnum\inputlineno>\LASTLINE%
    \ifnum\LASTLINE<0%
      \global\PAGETOP=\inputlineno%
    \fi%
    \global\LASTLINE=\inputlineno%
    \ifnum\INPSP=0%
      \ifnum\inputlineno>\PAGETOP%
      \fi%
    \else%
    \fi%
  \fi%
\fi}
\def\PUSH#1{%
\SRC%
\ifnum\INPSP=0 \global\let\INPSTACKA=\CurrentInput \else%
\ifnum\INPSP=1 \global\let\INPSTACKB=\CurrentInput \else%
\ifnum\INPSP=2 \global\let\INPSTACKC=\CurrentInput \else%
\ifnum\INPSP=3 \global\let\INPSTACKD=\CurrentInput \else%
\ifnum\INPSP=4 \global\let\INPSTACKE=\CurrentInput \else%
\ifnum\INPSP=5 \global\let\INPSTACKF=\CurrentInput \else%
               \global\let\INPSTACKX=\CurrentInput \fi\fi\fi\fi\fi\fi%
\gdef\CurrentInput{#1}%
\WinEdt{<+ \CurrentInput}%
\global\LASTLINE=0%
\ifSRCOK\fi%
\global\advance\INPSP by 1}
\def\POP{%
\ifnum\INPSP>0 \global\advance\INPSP by -1  \fi%
\ifnum\INPSP=0 \global\let\CurrentInput=\INPSTACKA \else%
\ifnum\INPSP=1 \global\let\CurrentInput=\INPSTACKB \else%
\ifnum\INPSP=2 \global\let\CurrentInput=\INPSTACKC \else%
\ifnum\INPSP=3 \global\let\CurrentInput=\INPSTACKD \else%
\ifnum\INPSP=4 \global\let\CurrentInput=\INPSTACKE \else%
\ifnum\INPSP=5 \global\let\CurrentInput=\INPSTACKF \else%
               \global\let\CurrentInput=\INPSTACKX \fi\fi\fi\fi\fi\fi%
\WinEdt{<-}%
\global\LASTLINE=\inputlineno%
\global\advance\LASTLINE by -1%
\SRC}
\def\INPUT#1{\relax}
\let\OldINCLUDE=\include
\def\include#1{%Always ".tex" file type!
\EJECT%
\PUSH{#1.tex}%
\OldINCLUDE{#1}%
\POP}
\let\originalxxxeverypar\everypar
\newtoks\everypar
\originalxxxeverypar{\the\everypar\expandafter\SRC}
\everymath\expandafter{\the\everymath\expandafter\SRC}
\let\zzzxxxbibliography=\bibliography
\def\bibliography#1{\PUSH{\jobname.bbl}\zzzxxxbibliography{#1}\POP}
\output\expandafter{\expandafter\SRCOKfalse\the\output}
\numberwithin{equation}{section}

\begin{document}

% \begin{center}
% \thispagestyle{empty}
% \begin{huge}
% Local Well-Posedness and Instability of Travelling Waves in a Chemotaxis Model\vspace{2cm}\\
% \end{huge}
% \begin{Large}
% \textbf{Martin Meyries}\vspace{0.5cm}\\Institut f\"ur Analysis, 
% Universit\"at Karls\-ruhe (TH)\\ 76128 Karlsruhe, Germany \vspace{1cm}\\ Tel.: +49 721 608 6519\\
% Fax.: +49 721 608 7650 \vspace{0.5cm} \\ Email: martin.meyries@math.uni-karlsruhe.de \vspace{1cm}\\Abbreviated Version of the Title: \vspace{0.5cm} \\ \textsl{Travelling Waves in a Chemotaxis Model} \vspace{3cm} \\Revised version from July 28, 2009
% 
% \end{Large}
% 
% \end{center}
% 
% \newpage
% \setcounter{page}{1}

\title{Local Well-Posedness and Instability of Travelling Waves in a Chemotaxis Model}

\author{Martin Meyries \vskip 10pt
Institut f\"ur Analysis,
Universit\"at Karls\-ruhe (TH) \vskip 0.1pt 76128 Karlsruhe, Germany
}

% \address{M. Meyries, Institut f\"ur Analysis,
% Universit\"at Karls\-ruhe (TH), 76128 Karlsruhe, Germany.}
 \email{martin.meyries@math.uni-karlsruhe.de}



%\date\today

\begin{abstract}
We consider the Keller-Segel model for chemotaxis with a nonlinear diffusion coefficent and a singular sensitivity function. We show the existence of travelling waves for wave speeds above a critical value, and establish local well-posedness in exponentially weighted spaces in a neighbourhood of a wave. A part of the essential spectrum of the linearization, which has unbounded coefficients on one half-axis, is determined. Generalizing the principle of linearized instability without spectral gap to fully nonlinear parabolic problems, we obtain nonlinear instability of the waves in certain cases.
\end{abstract}

\maketitle

%\centerline{\bf Preliminary version of \today}

%%%%%%%%%%%%%%%%%%%%%%%%%%%%%%%%%%%%%%%%%%%%%%%%%%%%%%%%%%%%%%%%%%%%%%%%%%

% \tableofcontents

\section{Introduction}
Chemotaxis denotes the directed movement of a cell species towards the gradient of a chemical. It is an important mechanism in biology and was discovered, for instance, in the context of stem cells and neurons. For the biological background on chemotaxis we refer to \cite{Eis04}. Chemotactic behaviour can be modelled by the (simplified) \textsl{Keller-Segel model}
\beq \label{eq:ksallg}
\begin{array}{rcl}
u_t & = & \left (k(u)u_x\right )_x - \left(u \phi(v)v_x\right )_x \\ 
v_t & = & \beta v_{xx} + b(u,v),
\end{array}
\eeq
where $x\in \R$. This system was first analyzed in \cite{KS71}. The scalar functions $u$ and $v$ denote the densities of the cell species and the chemical, respectively. The second summand of the first equation describes the chemotactic movement of the species towards the gradient of the chemical which is determined by the sensitivity function $\phi$ depending on $v$. In addition it is assumed that the species and the chemical diffuse nondegenerately, the species possibly nonlinearly, and that the chemical is produced or removed, as described by a cinetic term $b(u,v)$.

In mathematics, quasilinear, strongly coupled evolution equations of chemotaxis type (\ref{eq:ksallg}) have attracted a lot of attention over the last decades, as pattern formation and blow-up phenomena were discovered in such systems. We refer to the surveys \cite{Hor03, Hor04} and the references therein for an overview.

In this paper we investigate the existence and qualitative properties of travelling waves in (\ref{eq:ksallg}). We are concerned with a quite general  nonlinear diffusion coefficent $k$ and a singular sensitivity function $\phi(v) = \chi\frac{1}{v}$. Assuming $k$ to be bounded on bounded intervals, a singular sensitivity is necessary for the existence of travelling waves in (\ref{eq:ksallg}), see \cite{KS71, Sch05}. For linear diffusion, criteria for existence and nonexistence of travelling waves for such $\phi$ are given in \cite{NI91, Sch05}. Explicit travelling wave solutions are derived in \cite{HS04}. We are not aware of a treatment of the case of nonlinear diffusion in (\ref{eq:ksallg}) in a travelling wave context.

In \cite{NI91}, the waves are shown to be linearly unstable. Further, for a nonsingular sensitivity function $\phi$ and an additional growth term in the first equation of (\ref{eq:ksallg}), existence and linear (in)stability of travelling waves is proved in \cite{FMT00}. To our knowledge, there are no well-posedness and nonlinear (in)stability results in the literature for a singular sensitivity function.

Choosing a cinetic term $b$ similiar as \cite{Sch05}, we consider the system
\beq \label{eq:ks}
\begin{array}{rcl}
u_t & = &  \left (k(u)u_x\right )_x - \chi \left(u \frac{v_x}{v}\right )_x \\ 
v_t & = & v_{xx} + \gamma v - luv.
\end{array}
\eeq
It is no restriction to set here $\beta = 1$, since one can always rescale $x\ra x\sqrt{\beta}$ in space for $\beta>0$. Note that because of this, in relations for other coefficients absolute values will occur in the following. 

Throughout the paper, we make the following hypothesis on the coefficients.
\begin{itemize}
\item[\textbf{(H)}]\textsl{In system (\ref{eq:ks}) we have $k(u) = \alpha + d(u)$, where $\alpha>0$ is a constant and $d\in C^3(\R)$ is a nonnegative function. We further assume that $D(\cdot)=\int_0^\cdot \frac{d(\tau)}{\tau}\emph{\D} \tau$ is bounded at zero. The coefficients $\chi$, $\gamma$ and $l$ are strictly positive constants.}
\end{itemize}
Note that our condition on $d$ implies $d(0) = 0$. 

Searching for travelling waves, we work in a moving coordinate system $\xi = x-ct\in \R$ with constant speed $c>0$. Then with $'=\D/\D \xi$, the system (\ref{eq:ks}) transforms into
\beq \label{eq:eveq}
\begin{array}{rcl}
u_t & = & (k(u)u')' + cu' - \chi \left (u\frac{v'}{v} \right )' \\
v_t & = & v'' + cv' + \gamma v -luv.
\end{array}
\eeq 
A travelling wave is a zero of the right-hand side of (\ref{eq:eveq}) living in the space $C^2(\ol{\R})^2$ of functions converging at $\pm\infty$, see below for notation. For wave speeds $c> 2\sqrt{\gamma}$ we show the existence of a front-pulse wave $(u_*,v_*)$, i.e. $u_*(-\infty) >0$, $u_*(+\infty)= v_*(\pm\infty) = 0$ and $u_*,v_*>0$, by reducing the system to a single second order equation with a so-called KPP nonlinearity. Here we proceed similiar to \cite{NI91, Sch05}. 

For a nonlinear stability analysis, several severe problems arise due to the singular sensitivity function. 

Local well-posedness of (\ref{eq:eveq}) near a wave is not trivial at all, since the $v$-component of a wave converges exponentially fast to zero for $\xi\ra \pm\infty$ and thus the last term of the first equation becomes singular. We can only allow perturbations of the wave in exponentially weighted spaces, depending on the asymptotics of a front-pulse solution. Using \cite[Chapter 8]{Lun95}, for local well-posedness we show that the right-hand side of (\ref{eq:eveq}), considered as a map of perturbations of $(u_*,v_*)$, is $C^1$ with locally Lipschitz continuous derivative in an exponentially weighted space of continuous functions. We further show sectoriality of the linearization in each perturbation. For this, if $\chi$ is small compared to $\alpha$ we have to make an additional restriction on the lower bound of the wave speed $c$, but still obtain local well-posedness for wave speeds above a critical value. See (R1) below for details. 

The linearization of (\ref{eq:eveq}) is a coupled $2\times2$-system of second order ordinary differential operators whose coefficients are unbounded on the left half-axis. We determine a part of the essential spectrum of the linearization in a travelling wave by relating its Fredholm properties to the hyperbolicity of a corresponding first order constant coefficient matrix. The standard literature (e.g. \cite{San02}) assumes the coefficents to be bounded (in \cite{NI91}, the coefficients are bounded, in spite a pulse solution). Roughly speaking, we get rid of the unbounded coefficients by restricting the perturbations on the right half-axis. In the natural case that the minimal diffusion coefficient of the species $\alpha$ is less than or equal to 1, the diffusion coefficient of the chemical, we can choose the exponential weight on the right half-axis such that spectral values with positive real part occur. Otherwise we obtain such spectral values under certain restrictions on $\chi$ and $c$, see (R2) below.

By the above choice of the exponential weight, the wave and its translates will not be contained in the space of perturbations. Therefore our notion of nonlinear instability of a wave is the instability of a single equilibrium in the sense of Lyapunov (see Remark \ref{sec:reminst} for a detailed discussion). We show nonlinear instability of a wave in this sense by generalizing the principle of linearized instability without spectral gap on fully nonlinear evolution equations. For this purpose we show the applicability of \cite[Theorem 5.1.5]{Hen81} on abstract fully nonlinear equations, using optimal regularity results in weighted H\"older spaces from \cite{Lun95}. In a somewhat simpler context this can be found in \cite[Section 4.2]{Lun04}.

The paper is organized as follows. In Section 2 we describe the framework and state our main results. In Section 3 the proof for the existence of travelling waves is given, and in Section 4 we show local well-posedness around a wave. In Section 5 we determine a part of the spectrum of the linearized problem. Finally, in Section 6 we show the principle of linearized instability without spectral gap for abstract fully nonlinear equations.  \vspace{0.3cm}

\textsc{Notation.} Throughout the paper we denote generic positive constants by $C$. $\|\cdot\|$ always denotes the sup-norm. We set $\R_+=[0,+\infty)$, $\R_-= (-\infty,0]$. If for a function $f:\R\ra\R^N$ the limits $\lim_{\xi\ra\pm\infty}f(\xi)$ exist, we write $f(\pm\infty)$ for them. For functions $f,g:\R\ra \R$ we write $f\sim g$ as $\xi\ra \pm\infty$ if $f(\xi)/g(\xi) \ra 1$ as $\xi\ra \pm\infty$. If a function $f$ is $k$-times continuously differentiable, we often write simply $f\in C^k$. For $I=\R$,  resp. $\R_+$, we denote by $C^k(\ol{I})$ the set of functions $f\in C^k(I)$, where $f^{(j)}$ has finite limits at $\pm\infty$ and $+\infty$, respectively, for any integer $0\leq j\leq k$. Equipping $C^k(\ol{I})$ with the norm $\|f\|_{C^k} = \sum_{j=0}^k \|f^{(k)}\|$, it becomes a Banach space. If $X$ and $Y$ are Banach spaces, $B(X)$ and $B(X,Y)$ denote the Banach spaces of bounded linear operators on $X$ and from $X$ to $Y$, respectively. For $r>0$ and $x_0\in X$, $B_r^X(x_0)$ denotes the open ball of radius $r$ in $X$ with center $x_0$. \vspace{0.3cm}

\textsc{Acknowledgements.} This work is partially supported by Studienstiftung des deutschen Volkes, and partially included in my diploma thesis, supervised by Willi J\"ager at the University of Heidelberg. I am grateful to him for introducing me to the topic. I would also like to thank Bj\"orn Sandstede for helpful advices. Finally, I am deeply indebted to Roland Schnaubelt for encouragement and many helpful discussions.

\section{Framework and Results}\label{sectionfr}
For wave speeds above a critical value we obtain the existence of travelling waves in system (\ref{eq:ks}).
\begin{thm}\label{sec:ex}
\textsl{Assume (H). Then for each wave speed $c> 2\sqrt{\gamma}$ there exists a travelling wave solution $(u_*,v_*)$ of (\ref{eq:ks}). Its components $u_*,v_*$ are strictly positive, $u_*$ is shaped as a front, and $v_*$ is shaped as a pulse. More precisely, $u_*$ is strictly decaying, and $$u_*(-\infty) = u_*^-= l^{-1}(c^2/\chi^2 + c^2/\chi + \gamma)>0,\qquad  u_*(+\infty)= v_*(\pm\infty) = 0.$$}
\end{thm}
This theorem is proved in Section \ref{sectionex}, following the lines of \cite{Sch05}. Note that due to the scaling invariance of $v$ in (\ref{eq:ks}) and the translation invariance of the problem, for every wave speed $c> 2\sqrt{\gamma}$ Theorem \ref{sec:ex} gives a two parameter family of travelling wave solutions with positive components: 
\beq\label{eq:family}
\{(u_*,\lambda v_*)(\cdot+\xi_0)\;|\;\lambda>0, \;\xi_0\in \R\}.
\eeq
Fixing $c> 2\sqrt{\gamma}$ and a travelling wave $(u_*,v_*)$ as a zero of the right-hand side of (\ref{eq:eveq}), we make the following restriction on the wave speed.
\begin{itemize}
\item[\textbf{(R1)}] \textsl{In the case $\chi\leq \alpha-2$ it holds that}
\beq\label{eq:res1}
c^2>\gamma\frac{(\chi-\alpha)^2}{\alpha-\chi-1}.
\eeq
\end{itemize}
There is no restriction on $c$ if $\chi>\alpha-2$. 

Assuming (R1), $$J = \left [ \frac{c}{2} - \frac{1}{2} \sqrt{ c^2-4\gamma}, \; \frac{c}{\alpha} + \frac{c\chi}{2\alpha} - \frac{\chi}{2\alpha}\sqrt{c^2-4\gamma}\right ]$$ is an interval, containing more than one point. Now set $$w_-=-\frac{c}{\chi}, \qquad \trm{and take some }\quad   w_+\in J,$$ then $w_-<0$ and $w_+>0$. Note that for simplicity we will abbreviate $a= c/\chi$ in Section \ref{sectionex}, hence $w_-= -a$. Choose smooth functions $\eta_-, \eta_+\geq 1$ with the properties 
$$\eta_-(\xi) = \left \{ \begin{array}{ll} e^{w_-\xi}, & \xi<-1,\\ 1, & \xi\geq 0, \end{array}\right. \qquad \eta_+(\xi) = \left \{ \begin{array}{ll} 1, & \xi<0,\\ e^{w_+\xi}, & \xi\geq 1, \end{array}\right.$$ and set $\eta = \eta_-\cdot \eta_+$. Lemma \ref{sec:asym1} below shows that $\eta_-$ grows on $\R_-$ as $1/v_*$, and $\eta_+$ grows on $\R_+$ at least as $1/v_*$ and at most as $1/u_*$. Define $$X_1 = \{x\in C(\ol{\R})\,|\; \eta_+ x \in C(\ol{\R})\}, \qquad X_2 = \{y\in C(\ol{\R})\,|\; \eta y\in C(\ol{\R})\},$$ which are Banach spaces equipped with weighted norms $\|\cdot\|_{X_1} = \|\eta_+\cdot\|$ and $\|\cdot\|_{X_2} =\|\eta\cdot\|$, respectively, where $\|\cdot\|$ denotes, as always, the sup-norm. The space $X_1$ is not weighted on $\R_-$, corresponding to the fact that $u_*$ does not vanish for $\xi\ra - \infty$. We further define $$D_1 =  \{x\in C^2(\ol{\R})\,|\; x,x',x''\in X_1\},\qquad D_2 =  \{y\in C^2(\ol{\R})\,|\; y,y',y''\in X_2\}.$$
These are Banach spaces equipped with the weighted $C^2$-norms $\|\cdot\|_{D_1}$ and $\|\cdot\|_{D_2}$, respectively, where $\|x\|_{D_1} = \|x\|_{X_1} + \|x'\|_{X_1} + \|x''\|_{X_1}$, and analogously for $\|\cdot\|_{D_2}$. Finally, we set $$X= X_1\times X_2, \qquad D= D_1\times D_2,$$ equipped with the norms $\|\cdot\|_X = \|\cdot\|_{X_1} + \|\cdot\|_{X_2}$ and $\|\cdot\|_D = \|\cdot\|_{D_1} + \|\cdot\|_{D_2}$, respectively.
We consider the right-hand side of (\ref{eq:eveq}) as a map $F$ of perturbations $(x,y)$ of a travelling wave $(u_*,v_*)$. Being precise, we have 
\beq\label{eq:F}
F\left(\begin{array}{c} x\\ y\end{array} \right ) = \left(\begin{array}{c} (k(u_*+x)(u_*+x)')' + c(u_*+x)' - \chi \left ((u_*+x)\frac{(v_*+y)'}{v_*+y} \right)'\\ (v_*+y)'' + c(v_*+y)' + \gamma (v_*+y) -l(u_*+x)(v_*+y) \end{array} \right ).
\eeq 
Now $(u_*,v_*)$ corresponds to the zero $(0,0)$ of $F$. Proposition \ref{sec:F} shows that there is an open neighbourhood $\OO\subset D$ of $(0,0)$ such that $F:\OO\ra X$ is defined for any rate $w_+\in J$.

For an interval $[a,b]$, a Banach space $E$ and $\theta\in(0,1)$, we consider the H\"older space $C^\theta([a,b];E)$ equipped with the norm $\|\cdot\|_{C^\theta}= \|\cdot\|+ [\cdot]_{C^\theta}$, where $\|f\|= \sup_{t\in[a,b]}\|f(t)\|_E$ and $$[f]_{C^\theta([a,b]; E)} = \sup_{t,s\in [a,b], t>s} \frac{\|f(t)-f(s)\|_E}{(t-s)^\theta}.$$ For $\theta\in(0,1)$ we further introduce the weighted H\"older space
\beq\label{eq:Caa}
\begin{array}{rl}
C_\theta^\theta(]a,b]; E) = & \{f:]a,b]\ra E \trm{ bounded}\,| f\in C^\theta([a+\eps,b];E)\, \forall\, \eps\in (0,b-a),\\ 
&  \sup_{\eps\in(0,b-a)}\eps^\theta [f]_{C^\theta([a+\eps,b]; E)}<+\infty\},
\end{array}
\eeq
equipped with the norm $$\|f\|_{C_\theta^\theta(]a,b]; E)} = \|f\|+ \sup_{0<\eps<b-a}\eps^\theta [f]_{C^\theta([a+\eps,b]; E)},$$ cf. \cite[Chapter 4]{Lun95}. 

Combining \cite[Theorem 8.1.1 and Proposition 8.2.3]{Lun95}, we obtain the following abstract local well-posedness result. 
\begin{thm}\label{sec:LW}
\textsl{Assume that $X$ and $D$ are Banach spaces such that $D$ is continuously and densely embedded in $X$, that $\OO$ is an open neighbourhood of $\ue=0$ in $D$ and that $F:\mathcal{O}\ra X$ is a map with $F(0) = 0$, having the following properties:}
\begin{itemize}
 \item[(P1)] $F\in C^1(\OO,X)$,
 \item[(P2)] $F':\OO\ra B(D,X)$ \textsl{is locally Lipschitz continuous},
 \item[(P3)] $F'(\ue)$\textsl{, considered as an operator on $X$ with domain $D$, is sectorial for any $\ue\in \mathcal{O}$, and its graph norm is equivalent to the norm in $D$.}
\end{itemize}
\textsl{Then the evolution equation $\ue_t=F(\ue)$ is locally well-posed in $\mathcal{O}$. More precisely, for fixed $\theta\in (0,1)$ we have:} 
\begin{itemize}
\item[(LW1)] \textsl{For each $\ue_0\in \mathcal{O}$ there is an existence time $\tau(\ue_0)>0$ and a solution $$\ue(\cdot,\ue_0)\in C^1([0,\tau(\ue_0)];X)\cap C([0,\tau(\ue_0)]; D)\cap C_\theta^\theta(]0,\tau(\ue_0)]; D)$$ of $\ue_t(t) = F(\ue(t))$ for $t\in [0,\tau(\ue_0)]$ such that $\ue(0) = \ue_0$. The solution $\ue(\cdot,\ue_0)$ is unique in the set $$\bigcup_{0<\beta<1} C_\beta^\beta(]0,\tau(\ue_0)]; D)\cap C([0,\tau(\ue_0)]; D).$$
\item[(LW2)] For each given $T>0$ there is $\eps>0$, such that if $\ue_0\in B_\eps^D(0) \cap \OO$ then $\tau(\ue_0)\geq T$, and the solution map $$B_\eps^D(0) \cap \OO\ra C_\theta^\theta(]0,T]; D), \qquad \ue_0\mapsto \ue(\cdot,\ue_0),$$ is locally Lipschitz continuous.}
\end{itemize}
\end{thm}
By a sectorial operator we mean the generator of an analytic semigroup, see \cite[Chapter 2]{Lun95}.

In Section 4 we verify (P1)-(P3) for $F:\OO\ra X$ defined in (\ref{eq:F}) to obtain the following result.
\begin{thm}\label{sec:thmlw}
\textsl{Assuming (H) and (R1), the evolution equation
\beq\label{eq:eveq2}
\left(\begin{array}{l} x_t\\ y_t\end{array} \right ) = F\left(\begin{array}{l} x\\ y\end{array} \right )
\eeq
is locally well-posed for $(x,y)\in \OO$ in the sense (LW1), (LW2).}
\end{thm}

This shows that $F$ generates a dynamical system in $D$ in an open set of perturbations of a wave. Once this is established, we perform a stability analysis, considering the wave $(u_*,v_*)$ as the equilibrium $(0,0)$ of (\ref{eq:eveq2}). 

The linearization at the equilibrium, $F'(0,0)$, is a nondegenerate second order ordinary differental operator with continuous matrix coefficients, which are unbounded at $-\infty$ and converge at $+\infty$. In Section \ref{spectrum} we show that if $F'(0,0)$ is Fredholm considered as an operator on functions on the real line then it is Fredholm considered on functions on the right half-line. Now it is a standard result that the Fredholm properties of $F'(0,0)$ on $\R_+$ are closely related to the hyperbolicity of its corresponding first order constant coefficient operator. Thus we are able to calculate a part of the essential spectrum of $F'(0,0)$. Since we cannot determine the complete spectrum, we only make statements on the instability of the wave (see also Remark \ref{remtrunk}).

To obtain positive real parts in the spectrum of the linearization, we have to make the following assumptions.
\begin{itemize}
\item[\textbf{(R2)}] \textsl{It holds that $\chi>\alpha-2$. Further, if $\alpha>1$ then there is the upper bound 
\beq\label{eq:crest} 
c^2<\gamma \frac{(\alpha+\chi)^2}{(\chi+1)(\alpha-1)} 
\eeq  
on the wave speed $c$. In addition, the rate of the exponential weight on $\R_+$ is taken from $J_u$, where $J_u\subset J$ is the subinterval
\beq\label{eq:n12}
J_u = \left (\frac{c}{2}+ \frac{1}{2} \sqrt{c^2-4\gamma},\;  \frac{c}{\alpha} + \frac{c\chi}{2\alpha} - \frac{\chi}{2\alpha}\sqrt{c^2-4\gamma} \right].
\eeq}
\end{itemize}
One checks that the fraction in (\ref{eq:crest}) is always greater than $4$ for $\chi>\alpha-2$. In Theorem \ref{sec:ex}, the condition on the wave speed for the existence of travelling waves is $c^2> 4\gamma$. Thus for $\chi>\alpha-2$ there are always wave speeds such that (\ref{eq:crest}) holds. In applications, however, the species is expected to diffuse slower than the chemical, i.e. one assumes naturally that $\alpha\leq 1$.

In Section \ref{spectrum} we prove the following result for the spectrum.
\begin{thm}\label{sec:spec}
\textsl{Assuming (H) and (R2), it holds that $$\sigma(F'(0,0)) \cap \{\emph{\rea} \lambda >0\} \neq \emptyset.$$ More precisely, the curve $$\left \{\lambda \in \C\,|\,\emph{\rea} \lambda = -h^2 + w_+^2 - w_+c +\gamma,\;\emph{\ima} \lambda = (2w_+-c)h, \;h\in \R \right \},$$ which intersects the imaginary axis, is contained in $\sigma(F'(0,0)).$}
\end{thm}

\begin{rem}\label{sec:reminst} Due to Lemma \ref{sec:asym1} below, for the second component of the wave we have $v_*(\xi)\sim C\,e^{-(\min J)\xi}$ as $\xi\ra +\infty$, i.e. its exponential rate equals the lower bound for $w_+$ in $J$. Thus by assuming $w_+\in J_u$ in (R2), the wave and its translates are not contained in $D$, and the only (known) zero of $F$ in $\OO$ is $(0,0)$.
 
The most natural notion for nonlinear instability of a wave is orbital instability, i.e. the whole family (\ref{eq:family}) is unstable under perturbations. But this makes no sense in our setting when assuming that $w_+\in J_u$, since we cannot measure the distance of a perturbed wave to the translated waves in the $\eta_+$-weighted norm. Being precise, let $v_{**}(\cdot) =  v_*(\cdot+\xi_0)$ with $\xi_0\in \R$ be a translation of the second component of a wave and $y\in D_2$ be an arbitrary perturbation. Then, since $y$ decays faster than $v_*$ as $\xi\ra +\infty$ by assumption, $v_*+y \sim C\,e^{-(\min J)\xi}$ as $\xi\ra+\infty$. Therefore $$\eta_+(\xi)\left ( v_{**}(\xi)-v_{*}(\xi)-y(\xi)\right ) \sim e^{(w_+-(\min J))\xi}\left (C e^{a\xi_0} -  C \right )$$ does not converge for $\xi\ra +\infty$  for $\xi_0\neq 0$ if $w_+>\min J$.

\textsl{In this situation we define a wave to be nonlinearly unstable if it is nonlinearly unstable as a single equilibrium in $\OO$ in the sense of Lyapunov, see Theorem \ref{sec:henry} below or \cite[Section 2.9]{Per96} for a definition.}
\end{rem}

In Section \ref{sec:instability} we prove the principle of linearized instability without spectral gap for fully nonlinear parabolic problems.
\begin{thm} \label{sec:linin}
\textsl{In the setting of Theorem \ref{sec:LW}, suppose $F$ is in addition $p$-linearizable in $\ue=0$ for some $p>1$, and $\sigma(F'(0))\cap \{\emph{\rea}\lambda >0\}\neq \emptyset$.}

\textsl{Then the steady state $\ue=0$ of the locally well-posed evolution equation $\ue_t=F(\ue)$ is nonlinearly unstable in the sense of Lyapunov.}
\end{thm}
See (\ref{eq:plin}) for the definition of $p$-linearizability. For instance, this condition is fulfilled if $F\in C^2$. 

In Proposition \ref{eq:Fplin} we show that the right-hand side $F$ (see (\ref{eq:F})) is 2-linearizable. Thus we immediately obtain our final result.

\begin{thm}\label{sec:finthm}
\textsl{Assuming (H) and (R2), each travelling wave solution $(u_*,v_*)$ from Theorem \ref{sec:ex} is nonlinearly unstable in the sense of Lyapunov in the exponential weighted space $D$ with respect to perturbations $(x,y)\in \OO$.}
\end{thm}

\section{Proof of Theorem \ref{sec:ex}: Existence of Travelling Waves}\label{sectionex}
Denoting the wave speed by $c>0$, we are searching for nonnegative solutions which are constant in the moving frame $\xi = x-ct\in \R$, i.e. nonnegative functions $u_*,v_*\in C^2(\ol{\R})$ such that  $u(x,t) = u_*(x-ct)$ and $v(x,t)=v_*(x-ct)$ solve (\ref{eq:ks}) for $x\in \R$ and $t\in \R$.
% \begin{thm}\label{sec:ex}
% \textsl{Assume (H). Then for every wave speed $c\geq 2\sqrt{\gamma}$, there exists a travelling wave solution $(u_*,v_*)$ of (\ref{eq:ks}). Its components $u_*,v_*$ are strictly positive, $u_*$ is shaped as a front, and $v_*$ is shaped as a pulse. More precisely, $u_*$ is strictly decaying and $$u_*(-\infty) = u_*^-= l^{-1}(c^2/\chi^2 + c^2/\chi + \gamma)>0,\quad  u_*(+\infty)= 0,\quad v_*(\pm\infty) = 0.$$}
% \end{thm}

Our strategy follows \cite{NI91, Sch05}. Plugging the travelling wave ansatz into (\ref{eq:ks}) and writing $'=\D/\D \xi$, we obtain the following system of ordinary differential equations:
\begin{eqnarray} \label{eq:twansatz1}
-cu' & = & \left ((\alpha + d(u))u' - \chi u\frac{v'}{v}\right )'\\ \label{eq:twansatz2}
-cv' & = & v'' + \gamma v - luv
\end{eqnarray}
First we manipulate this system rather informal, collecting information needed for a rigorous existence proof. 

Supposing $v>0$, integrating (\ref{eq:twansatz1}) and neglecting constants of integration yields $$u'= \frac{u}{\alpha+d(u)}( \chi (\log v)' - c).$$ Setting $D(\cdot)=\int_0^\cdot \frac{d(\tau)}{\tau}\D \tau$ and supposing $u>0$, we obtain 
\begin{equation}\label{eq:n1}
\alpha\log u(\xi) + D(u(\xi)) = \chi \log v(\xi) - c\xi.
\end{equation}

\begin{lem}
\textsl{The map $G:(0,+\infty) \ra \R$, $G(u) = \alpha \log u + D(u)$, is bijective, strictly increasing, and $C^2$. Its inverse $G^{-1}$ is $C^2$, satisfies
\begin{equation}\label{eq:Gasym}
G^{-1}(y) \sim e^{y/\alpha}\qquad \trm{as }\; y\ra -\infty
\end{equation}
and, for $y\in \R$, 
\begin{equation}\label{eq:ablGinv}
(G^{-1})'(y) = \frac{G^{-1}(y)}{\alpha + d(G^{-1}(y))}.
\end{equation}}
\end{lem}
\bprf
We calculate $G'(u) = (\alpha+d(u))/u$ and see that $G\in C^2$. From $\alpha+d(\cdot)>0$ it follows that $G'>0$, therefore $G$ is injective. From $D(0) = 0$ and $D\geq 0$  we conclude that $G(u)\ra -\infty$ for $u\ra 0$ and $G(u) \ra +\infty$ for $u\ra +\infty$. This shows the invertibility of $G$, so its inverse $G^{-1}$ exists. For any $\eps>0$ there is a $C>0$, such that $e^{y/\alpha-\eps}\leq G^{-1}(y) \leq e^{y/\alpha}$ for $y<-C$, so (\ref{eq:Gasym}) follows. Since $G'$ never vanishes, the inverse is $C^2$, with the stated formula for $(G^{-1})'$.
\eprf Now we can solve in (\ref{eq:n1}) for $u$ in terms of $v$ and $\xi$, obtaining $u(\xi) = G^{-1}(\chi\log v(\xi) - c\xi).$ Plugging this into (\ref{eq:twansatz2}) yields 
\beq\label{eq:star}
v'' + cv' + v(\gamma - lG^{-1}(\chi\log v - c\xi)) = 0.
\eeq We abbreviate 
\beq \label{eq:munu}
a=c/\chi, \qquad \mu = 2a +c\qquad \nu = a^2 + ac+\gamma.
\eeq Then by setting $v(\xi) = e^{a\xi}p(\xi)$ with an unknown positive scalar function $p$, equation (\ref{eq:star}) is transformed to
\beq \label{eq:n2}
p'' + \mu p' + f(p) = 0,
\eeq 
where $f:\R\ra \R$ is defined by 
\beq \label{eq:eff}
f(p) = \left \{  \begin{array}{ll} p(\nu-lG^{-1}(\chi \log|p|)), & p\neq 0, \\ 0, & p=0. \end{array} \right.
\eeq The next lemma states that $f$ is a so-called KPP nonlinearity (cf. \cite{KPP37}).
\begin{lem}\label{sec:f}
\textsl{The function $f$ is $C^1$ on $\R$ and has the properties $$f(0) = f(p_0) = 0,\qquad  f(p)>0 \quad \trm{for } p\in (0,p_0), \quad f'(0)=\nu >0, \quad f'(p_0) <0,$$ where $p_0 = \exp\left( \frac{G(\nu/l)}{\chi}\right)>0.$}
\end{lem}
\bprf Clearly $f$ vanishes at the points $0$ and $p_0$. Since $G^{-1}(-\infty) = 0$ and $G^{-1}$ grows strictly, $f$ is positive on the interval $(0,p_0)$. We have $f \in C^1$ for $p\neq 0$ with derivative $$f'(p) = \nu - lG^{-1}(\chi \log|p|)) -\chi l(G^{-1})'(\chi \log|p|).$$ From (\ref{eq:ablGinv}) we deduce that $f'(p)\ra \nu>0$ for $p\ra 0$. Further $f'(0) = \nu$, so $f\in C^1$. Plugging $p_0$ into $f'$ yields $f'(p_0)<0$.
\eprf
We rewrite (\ref{eq:n2}) as the first order system
\begin{equation} \label{eq:autoeo}
\left( \begin{array}{cc} p'\\ q' \end{array} \right ) = \left( \begin{array}{c} q \\  -f(p) -\mu q \end{array} \right ).
\end{equation} 
The system possesses the steady states $(p_0,0)$ and $(0,0)$. For certain values of $\mu$, the following classical result states the existence of a heteroclinic orbit of (\ref{eq:autoeo}) connecting the two steady states.
\begin{thm}[{\cite[IV.2.3]{DT76}}]\label{sec:DT76} \textsl{There exists a $\mu_0> 0$ with 
\beq\label{eq:muecond} 
2\sqrt{f'(0)} \leq \mu_0 \leq 2\sqrt{\sup_{0<p<p_0}\frac{f(p)}{p}},
\eeq such that (\ref{eq:autoeo}) possesses a heteroclinic orbit $(p_*,q_*)$ with $(p_*,q_*)(-\infty) = (p_0,0)$,  $(p_*,q_*)(+\infty) = (0,0)$, $p_*>0$ and $q_*<0$, provided that $\mu> \mu_0$.}
\end{thm}
The special form of $f$ and (\ref{eq:muecond}) imply that $\mu_0=2\sqrt{\nu}$. In what follows, let $p_*$ and $q_*$ be given by Theorem \ref{sec:DT76}, where we assume that $\mu>\mu_0$. This last condition is equivalent to $c> 2\sqrt{\gamma}$. The linearization $H$ of the right-hand side of (\ref{eq:autoeo}) in $(0,0)$, $$H=\left( \begin{array}{cc} 0 & 1 \\ -\nu & -\mu \end{array} \right ),$$
is hyperbolic with differing real eigenvalues 
\beq \label{eq:ev}
\kappa_\pm = \frac{-\mu \pm \sqrt{\mu^2-4\nu}}{2} = -a -\frac{c}{2} \pm \frac{1}{2}\sqrt{c^2-4\gamma}<0.
\eeq
We investigate the asymptotic behaviour of the first component of the heteroclinic orbit for $\xi\ra +\infty$ in detail. This is not really needed in the existence proof, but is crucial in the next sections. We rewrite (\ref{eq:autoeo}) as
\beq \label{eq:linsys}
\left ( \begin{array}{c} p'\\ q'\end{array} \right )  = H \left ( \begin{array}{c} p\\ q\end{array} \right ) + \left ( \begin{array}{c} 0 \\ g(p)\end{array} \right ),
\eeq
where the function $g:\R\ra\R$ defined by 
$$
g(p)=\left \{  \begin{array}{ll} lp G^{-1}(\chi\log|p|), & p\neq 0, \\ 0, & p=0. \end{array} \right.
$$
Then $g$ is $C^1$ and has the properties $g(\R_+)\subset \R_+,$ $g(0) = 0,$ $g(p_0)=p_0\nu >0$, $g'(0)=0,$ as one verifies as in the proof of Lemma \ref{sec:f}.
\begin{lem}\label{sec:key1}
\textsl{There is a constant $S_1>0$, such that $e^{-\kappa_+\xi}p_*(\xi) \ra S_1$ as $\xi\ra +\infty$. Further, $e^{-\kappa_+ \xi} q_*(\xi)$ converges for $\xi\ra +\infty$, and $(q_*/p_*)(+\infty) = \kappa_+$.}
\end{lem}
\bprf
For simplicity we neglect the stars for the heteroclinic orbit and write $(p,q)$. Using (\ref{eq:linsys}), the orbit can implicitly be represented by the variation of constants formula (\cite[Section 1.10]{Per96}). The eigenvectors of $H$ corresponding to its eigenvalues $\kappa_\pm$ are $(1,\kappa_\pm)$. We diagonalize $H$ by $H= TDT^{-1}$, where $$T = \left ( \begin{array}{cc} 1 & 1 \\ \kappa_+ & \kappa_- \end{array} \right ), \quad T^{-1} = \frac{1}{\kappa_- - \kappa_+}\,\left  ( \begin{array}{cc} 
\kappa_- & -1 \\ -\kappa_+ & 1 \end{array} \right ), \quad D= \left  ( \begin{array}{cc} 
\kappa_+ & 0 \\ 0 & \kappa_- \end{array} \right ),$$ and obtain $$e^{H\xi} = T e^{D\xi} T^{-1} = \frac{1}{\kappa_- - \kappa_+} \left ( \begin{array}{cc} \kappa_- e^{\kappa_+\xi} - \kappa_+ e^{\kappa_-\xi} & -e^{\kappa_+\xi}+ e^{\kappa_-\xi} \\ \kappa_+\kappa_- e^{\kappa_+\xi} -\kappa_+\kappa_- e^{\kappa_-\xi} & -\kappa_+ e^{\kappa_+\xi} + \kappa_- e^{\kappa_-\xi} \end{array} \right ).$$ Choosing an arbitrary point $(p(0),q(0))$ on the heteroclinic orbit $(p,q)$, the variation of constants formula yields
\begin{align*}
\left ( \begin{array}{c} p(\xi) \\ q(\xi) \end{array} \right ) = & \frac{1}{\kappa_- - \kappa_+} \left ( \begin{array}{cc} \kappa_- e^{\kappa_+\xi} - \kappa_+ e^{\kappa_-\xi} & -e^{\kappa_+\xi}+ e^{\kappa_-\xi} \\ \kappa_+\kappa_- e^{\kappa_+\xi} -\kappa_+\kappa_- e^{\kappa_-\xi} & -\kappa_+ e^{\kappa_+\xi} + \kappa_- e^{\kappa_-\xi} \end{array} \right ) \\
& \cdot \left ( \begin{array}{c} p(0) + \int_0^\xi \frac{1}{\kappa_- - \kappa_+} \left ( -e^{-\kappa_+s}+ e^{-\kappa_-s}\right) g(p(s))\D s  \\ q(0) + \int_0^\xi \frac{1}{\kappa_- - \kappa_+}\left (  -\kappa_+ e^{-\kappa_+s} + \kappa_- e^{-\kappa_-s} \right ) g(p(s))\D s \end{array} \right ).
\end{align*} 
After a careful calculation we obtain
\begin{align*}
(\kappa_- - \kappa_+) p(\xi) =  & \left (  \kappa_- e^{\kappa_+\xi} - \kappa_+ e^{\kappa_- \xi} \right ) p(0)  & + &  \left ( -e^{\kappa_+\xi} + e^{\kappa_-\xi}\right )q(0)\\
 &  -  \int_0^{\xi} e^{\kappa_+(\xi-s)}g(p(s))\D s  & + & \int_0^{\xi} e^{\kappa_-(\xi-s)}g(p(s))\D s, \\
 (\kappa_- - \kappa_+) q(\xi) = & \; \kappa_+\kappa_-\left (   e^{\kappa_+\xi} -e^{\kappa_- \xi} \right ) p(0) & + &  \left ( -\kappa_+e^{\kappa_+\xi} +\kappa_- e^{\kappa_-\xi}\right )q(0)\\
  &  -  \kappa_+\int_0^{\xi} e^{\kappa_+(\xi-s)}g(p(s))\D s & + & \;\kappa_-\int_0^{\xi} e^{\kappa_-(\xi-s)}g(p(s))\D s 
\end{align*}
for the components. Since $(p,q)$ is an orbit on the stable manifold of the hyperbolic steady state $(0,0)$, its components approach $(0,0)$ with an exponential rate larger, but arbitrarily close to $\kappa_+<0$ as $\xi\ra +\infty$ (\cite[Section 2.7]{Per96}). Using (\ref{eq:Gasym}) and the formula for $g$, this implies that $g(p(s))\leq C\,e^{(1+\frac{\chi}{\alpha})(\kappa_+ +\eps)s}$ for arbitrary small $\eps>0$ if $s$ is chosen large enough. Using $\frac{\chi}{\alpha} >0$ and multiplying the formulas for $p$ and $q$ by $e^{-\kappa_+\xi}$, we obtain the convergence of $e^{-\kappa_+\xi}p(\xi)$ and $e^{-\kappa_+\xi}q(\xi)$ to nonnegative values as $\xi\ra +\infty$, respectively. We now show that $e^{-\kappa_+\xi}p(\xi)$ is strictly growing, thus obtaining a limit $S_1>0$ as $\xi\ra +\infty$. Since $\left ( e^{-\kappa_+\xi}p(\xi)\right )' = e^{-\kappa_+\xi}(-\kappa_+ p(\xi) + q(\xi))$, we have to show that $q> \kappa_+ p$. Using the implicit representations of $p$ and $q$, we find after a careful calculation that this inequality is equivalent to 
\beq\label{eq:ineqpq}
\kappa_+p(0) - q(0) < \int_0^\xi e^{-\kappa_-s}g(p(s))\D s.
\eeq
Since $\kappa_+<0$, $q(-\infty)=0$ and $p(-\infty)=p_0>0$, the left-hand side of (\ref{eq:ineqpq}) will be strictly negative if we choose $(p(0),q(0))$ sufficiently close to the steady state $(p_0,0)$. The right-hand side is strictly positive for $\xi>0$. Hence $q>\kappa_+p$, and the first assertion is proved. For the last assertion, we use the product rule to obtain $$\frac{q(\xi)}{p(\xi)} = \kappa_+ + \frac{\left(e^{-\kappa_+\xi}p(\xi)\right)'}{e^{-\kappa_+ \xi}p(\xi)}.$$ The denominator on the right-hand side converges to $S_1>0$ as $\xi\ra +\infty$. Using (\ref{eq:eff}) and (\ref{eq:autoeo}), we see that 
$$\left(e^{-\kappa_+ \xi}p(\xi)\right)''=e^{-\kappa_+\xi}\left (\kappa_+^2p(\xi) -f(p) - (\mu+2\kappa_+)q(\xi)\right)$$
 is bounded on $\R_+$, hence $\left(e^{-\kappa_+\xi}p(\xi)\right)'$ is uniformly continuous on $\R_+$, and therefore $\left(e^{-\kappa_+\xi}p(\xi)\right)'\ra 0$ as $\xi\ra+\infty$.
\eprf
We are now ready to prove the existence result.
\begin{proof}[\emph{\textbf{Proof of Theorem \ref{sec:ex}.}}] System (\ref{eq:autoeo}) possesses a heteroclinic orbit $(p_*,q_*)$ as stated in Theorem \ref{sec:DT76} if and only if $\mu> 2\sqrt{\nu}$, which reads $c> 2\sqrt{\gamma}$. We claim that 
\beq\label{eq:uv}
u_*(\xi) = G^{-1}(\chi\log p_*(\xi)), \qquad v_*(\xi) = e^{a\xi} p_*(\xi)
\eeq 
are the components of a travelling wave solution as stated in Theorem \ref{sec:ex}. Using (\ref{eq:ablGinv}), (\ref{eq:munu}) and (\ref{eq:n2}) one easiliy verifies that $(u_*(\xi),v_*(\xi))$ solves (\ref{eq:twansatz1}), (\ref{eq:twansatz2}). Since $G^{-1}$ and $p_*$ are strictly positive, both components are strictly positive. We have 
\beq\label{eq:upr}
u_*' = \frac{G^{-1}(\chi\log p_*)}{\alpha + d(G^{-1}(\chi\log p_*))}\cdot \frac{\chi q_*}{p_*},
\eeq 
and since $G^{-1}>0$, $q_*<0$ and $p_*>0$, we conclude that $u_*$ decays strictly. It remains to check the asymptotic properties of the components. Since $p_*(-\infty)=p_0$, we conlcude that $v_*(-\infty)=0$ and that $u_*(-\infty)$ is as stated in the theorem. Lemma \ref{sec:key1} shows that $p_*(\xi)\sim S_1e^{\kappa_+\xi}$ for $\xi\ra+\infty$. Since $\kappa_+<-a$, see (\ref{eq:ev}), we conclude that $v_*(+\infty) = 0$. Finally, $u_*(+\infty) = 0$ follows from $G^{-1}(-\infty) =0.$
\eprf
For later purposes, we investigate the asymptotic behaviour of $u_*$, $v_*$ and their derivatives in more detail.
\begin{lem}\label{sec:asym1}
\textsl{There are constants $S_2,S_3>0$ such that the following holds. \\As $\xi\ra -\infty$:
$$
\begin{array}{lll}
u_*(-\infty) >0, & u_*'(-\infty) = 0, & u_*''(-\infty) = 0;
\end{array}
$$
as $\xi\ra +\infty$:
$$
\begin{array}{lll}
u_*(\xi) \sim S_2 e^{\left(\frac{\chi}{\alpha}\kappa_+\right)  \xi}, & u_*'(\xi)\sim S_3 e^{\left(\frac{\chi}{\alpha}\kappa_+\right)  \xi}, & e^{-\left(\frac{\chi}{\alpha}\kappa_+\right)  \xi}u_*''(\xi) \trm{ converges;}
\end{array}
$$
as $\xi\ra -\infty$:
$$
\begin{array}{lll}
v_*(\xi) \sim p_0e^{a\xi}, & e^{-a\xi}v_*'(\xi) \trm{ converges,} & e^{-a\xi}v_*''(\xi) \trm{ converges;}
\end{array}
$$
as $\xi\ra +\infty$:
$$
\begin{array}{lll}
v_*(\xi) \sim S_1 e^{(\kappa_++a)\xi}, & e^{-(\kappa_++a)\xi}v_*'(\xi)\trm{ converges,} & e^{-(\kappa_++a)\xi}v_*''(\xi) \trm{ converges.}
\end{array}
$$
Note that $a$ and $\kappa_+$ are defined in (\ref{eq:munu}) and (\ref{eq:ev}).}
\end{lem}
\bprf
Recall (\ref{eq:uv}). Setting $y= \chi\log( p_*(\xi))$ in (\ref{eq:Gasym}) and using Lemma \ref{sec:key1}, we obtain the asymptotics for $u_*$ as $\xi\ra+\infty$. The derivative $u_*'$ was calculated in (\ref{eq:upr}). We see that $u_*'(-\infty) = 0$. Using $(q_*/p_*)(+\infty)=\kappa_+$, we obtain the asserted asymptotics for $u_*'$ as $\xi\ra+\infty$. Differentiating (\ref{eq:upr}) once more, we obtain $$u_*'' = \frac{\chi q_*}{p_*} \left (\frac{u_*'}{\alpha+d(u_*)} - \frac{u_*u_*'d'(u_*)}{(\alpha+d(u_*))^2}\right ) + \frac{\chi u_*}{\alpha+d(u_*)}\left ( \frac{q_*'}{p_*} - \left(\frac{q^{*}}{p^{*}}\right)^2\right).$$ Using that $q_*' = - f(p_*)-\mu q_*$ and the Lemmas \ref{sec:f} and \ref{sec:key1}, we can verify the assertions for $u_*''$. The claims for $v_*$, $v_*'$ and $v_*''$ easily follow from (\ref{eq:eff}), (\ref{eq:uv})  and Lemma \ref{sec:key1}.
\eprf

\section{Proof of Theorem \ref{sec:thmlw}: Local Well-Posedness}\label{sectionwp}
In this section we are forced to assume $-\frac{\chi}{\alpha}\kappa_+ + (a+\kappa_+)>0$, see the proof of Lemma \ref{sec:asym2} below. Performing elementary manipulations using (\ref{eq:munu}) and (\ref{eq:ev}), we see that this always holds if $\chi\geq \alpha-2$, and otherwise we have to assume (R1), formulated in Section \ref{sectionfr}. Using the abbrevations introduced in the last section, we obtain for the range of the exponential rate $w_+$ on $\R_+$
\beq\label{eq:J}
J = \left[-(a+\kappa_+), -\frac{\chi}{\alpha}\kappa_+\right].
\eeq
Assuming (R1), $J$ really is an interval, containing more than one point.

We collect some properties of the weighted space $D$. Note that, as for any exponential weight, $x\in D_1$ is equivalent to $\eta_+ x \in C^2(\ol{\R})$. The analogous property holds for $y \in D_2$.  It is easy to see that $D$ is continuously and densely embedded in $X$.

\begin{lem}\label{sec:U2}
\textsl{The set $\{y\in D_2\;|\; v_*+y >0\}$ contains an open neighbourhood $U_2$ of $y=0$.} \textsl{For each $y\in U_2$ we have $\|y\|_{D_2}\,<\min\{p_0/2,S_1/2\}$. Further, there is an $\eps>0$ such that $\eta(v_*+y)\geq \eps$ for each $y\in U_2$.}
\end{lem}
\bprf
Due to Lemma \ref{sec:asym1}, for any choice of $w_+$ there is a number $\xi_0>0$ such that $\eta v_* > m_1 =\min\{p_0/2,S_1/2\}$ for $|\xi|>\xi_0$. Set $m_2= \inf_{[-\xi_0,\xi_0]} v_*>0$ and $\eps = \frac{1}{2}\min\{m_1,m_2\}$. Hence, $\eta v_* + \eta y > \eta v_* - \eps\geq \eps$ on $\R$ for each $y\in U_2:= B_{\eps}^{D_2}(0)$.
\eprf
\begin{lem} \label{sec:fraccon}
 \textsl{For any $y\in U_2\subset D_2$ it holds that $$\frac{v_*' + y'}{v_*+y}, \qquad \frac{v_*'' + y''}{v_*+y}\qquad \in \quad C(\ol{\R}).$$}
\end{lem}
\bprf
The first fraction is equal to $$\frac{e^{a\xi}(ap_*(\xi) + q_*(\xi)) + y'(\xi)}{e^{a\xi}p_*(\xi) + y(\xi)}.$$ Expanded by $\eta_-(\xi)=e^{-a\xi}$, all occuring terms converge as $\xi\ra -\infty$. Since $p_*(-\infty)=p_0$ and $\|\eta_- y\|\leq \|y\|_{D_2}<p_0/2$, the denominator does not converge to zero.  Expanded by $e^{-(a+\kappa_+)\xi}$, as $\xi\ra+\infty$  all occuring terms converge by Lemma \ref{sec:key1}. The denominator does not converge to zero, since $e^{-\kappa_+\xi}p_*(\xi)\ra S_1$ and $|e^{-(a+\kappa_+)\xi}y(\xi)|\leq |\eta_+(\xi) y(\xi)|\leq S_1/2$ for $\xi\ra +\infty$. The second fraction is equal to $$\frac{e^{a\xi}(a^2p_*(\xi) + (2a-\mu)q_*(\xi) - f(p_*(\xi))) + y''(\xi)}{e^{a\xi}p_*(\xi) + y(\xi)}.$$ Recall that $f(p) = p(\nu-lG^{-1}(\chi\log p))$ for $p>0$. As before, expanded by $\eta_-$ one deduces the convergence as $\xi\ra -\infty$. Finally, expanded by $e^{-(a+\kappa_+)\xi}$, as $\xi\ra +\infty$ the convergence follows again from Lemma \ref{sec:key1}.
\eprf
With the help of Lemma \ref{sec:U2} we define the open neighbourhood $$\OO= D_1\times U_2$$ of $(0,0)$ in $D$. 

\begin{prop}\label{sec:F}
\textsl{The map $F:\OO\ra X$, where $$F\left(\begin{array}{c} x\\ y\end{array} \right ) = \left(\begin{array}{c} (k(u_*+x)(u_*+x)')' + c(u_*+x)' - \chi \left ((u_*+x)\frac{(v_*+y)'}{v_*+y} \right)'\\ (v_*+y)'' + c(v_*+y)' + \gamma (v_*+y) -l(u_*+x)(v_*+y) \end{array} \right )$$ as in (\ref{eq:F}), is defined.}
\end{prop}
\bprf
Recall Lemma \ref{sec:asym1}. Let $(x,y)\in \OO$. In the first component of $F(x,y)$, the first two summands are in $X_1$ since $k(0)=\alpha$ and the exponential rate of $u_*$ and its derivatives at $+\infty$ is larger or equal to $w_+$ for any choice of $w_+\in J$, see (\ref{eq:J}). The third summand is equal to $$-\chi(u_*'+x')\,\frac{v_*'+ y'}{v+y} -\chi (u_*+x) \frac{v_*''+ y''}{v+y} +\chi (u_*+x) \left ( \frac{v_*'+ y'}{v+y}\right )^2.$$ Due to Lemma \ref{sec:fraccon}, the fractions have limits at $\pm\infty$. The other factors are in $X_1$ as explained before. Using that $v_*$ satisfies the travelling wave equation (\ref{eq:twansatz2}), the second component of $F(x,y)$ is equal to $$y''+ c y'- lv_*x + y(\gamma-l(u_*+x)).$$ The third summand is in $X_2$ since $\eta_-v_*$ converges as $\xi\ra -\infty$ and $\eta_+x$ converges as $\xi\ra +\infty$. The other summands belong to  $X_2$ since $y\in D_2$ and $(u_*+x)\in C(\ol{\R})$.
\eprf

% \textsl{there is a solution map $$\Phi:\{(t,u_0)\,|\, u_0\in \OO, t\in [0,\tau(u_0)]\}\ra W(0,\tau(,$$ such that for every $u_0\in \OO$ $$\Phi (\cdot,u_0)\in$$ $u(\cdot): = \Phi(\cdot,u_0)$ is in  $C^1([0,\tau(u_0)];X)$ with $u(0) = u_0$ and $u(\cdot)$ is unique solution of $u_t(t) = F(u(t))$ for $t\in [0,\tau(u_0)]$;

% \begin{itemize}
%  \item[(LW1)] for any $u_0\in \OO$ there is an existence time $\tau(u_0)>0$ and $u\in C^1([0,\tau(u_0)];X)\cap C([0,\tau(u_0)]; D)$ such that $u_t(t) = F(u(t))$ for any $t\in [0,\tau(u_0)]$ and $u(0) = u_0$;
%  \item[(LW2)] for any given $T>0$ there exists $\eps>0$ such that if $u_0\in \OO$ with $\|u_0\|_D\leq \eps$, then $\tau(u_0)\geq T$ and 
% \beq\label{eq:loclip}
% \|u(\cdot,u_0)\|_{C_\alpha^\alpha(]0,T]; D)}\leq H\;\|u_0\|_D,
% \eeq where $H$ is a positive constant only depending on $T$.
% \end{itemize}}

To show local well-posedness of $(x,y)_t= F(x,y)$ in the sense (LW1)-(LW2) in $\OO$ for $(x,y)\in \OO$, we verify the assumptions on $X$, $D$ and $\OO$ and the properties (P1)-(P3) for $F$ from Theorem \ref{sec:LW}. 

It is clear that $D\subset X$ is densely embedded, $\OO\subset D$ is open, $(0,0)\in \OO$ and $F(0,0)=(0,0)$. Formally linearizing $F$ at $(x,y)\in \OO$, we obtain the operator
\begin{equation} \label{eq:A}
L_{(x,y)}=\left (\begin{array}{cc} L_1 & L_2 \\ L_3 & L_4 \end{array} \right ), 
\end{equation}
where 
\begin{align*}
&L_1= a_1\partial_{\xi\xi} +a_2 \partial_\xi + a_3,\quad L_2 = a_4\partial_{\xi\xi} + a_5 \partial_\xi + a_6,\\
&L_3 = a_7, \qquad \qquad   L_4 =  \partial_{\xi\xi} + c\partial_\xi + a_8,
\end{align*}
with the coefficents (we write $(u,v) = (u_*,v_*) + (x,y)$)
\begin{align}
a_1 &=\alpha+d(u) ,\quad a_2= c+2u'd'(u)-\chi\frac{v'}{v}, \nonumber \\
a_3&= d'(u)u'' +d''(u)(u')^2 -\chi \frac{v''}{v} -\chi v'\left(\frac{1}{v}\right)',\qquad a_4=-\chi\frac{u}{v},\label{eq:coeff} \\
a_5 &= -\chi \left (\frac{u'}{v} +  2u\left(\frac{1}{v}\right)'\right), \qquad a_6= -\chi\left(u'\left(\frac{1}{v}\right)' + u\left(\frac{1}{v}\right)''\right), \nonumber\\
a_7& = - lv, \qquad a_8 = \gamma - lu.\nonumber
\end{align}
We investigate the asymptotic properties of these coefficents and write $a_i^\pm=a_i(\pm\infty)$. 
\begin{lem}\label{sec:asym2}\textsl{Assume (R1), i.e. $\frac{\chi}{\alpha}\kappa_+ - (a+\kappa_+)<0$. Then for any $w_+\in J$ and arbitrary $(x,y)\in \OO$ the limits $a_1^\pm, a_2^\pm, a_3^\pm, a_7^\pm$, $a_8^\pm$ exist in $\R$, where $a_7(\xi)\sim -lp_0e^{a\xi}$ for $\xi\ra -\infty$. Further, $a_4^+$, $a_5^+$, $a_6^+$ exist but $a_4$, $a_5$, $a_6$ are unbounded on $\R_-$ and  $$e^{a\xi} a_i(\xi)\quad \trm{ converges as }\quad \xi\ra -\infty,  \qquad i=4,5,6.$$ In the case $(x,y)=(0,0)$ we obtain 
\beq\label{eq:apmlimits}
a_1^+ = \alpha,\quad a_2^+ = c + \frac{\chi}{2}\left (c - \sqrt{c^2-4 \gamma}\right),\quad a_3^+ = ... = a_7^+=0, \quad a_8^+ = \gamma.
\eeq} 
\end{lem}
\bprf The assertions for $a_1, a_2, a_3, a_7$ and $a_8$ follow from Lemma \ref{sec:asym1} and Lemma \ref{sec:fraccon}. The term  $$e^{a\xi}\,\frac{u(\xi)}{v(\xi)} = \frac{u_*(\xi)+x(\xi)}{p_*(\xi)+e^{-a\xi}y(\xi)}$$ converges for $\xi\ra-\infty$, since the denominator converges to a nonzero number due to  the choice of $y\in U_2$. Since $u_*'$ is given by (\ref{eq:upr}), we conclude the convergence of $e^{a\xi}u'(\xi)/v(\xi)$ for $\xi\ra -\infty$ in an analogous way. Calculating $(1/v)'$ and $(1/v)''$ we see that all assertions for $\xi\ra -\infty$ follow, since $a_4$, $a_5$ and $a_6$ consist of sums and products of $u/v$, $u'/v$ and fractions considered in Lemma \ref{sec:fraccon}. Next consider $$\frac{u(\xi)}{v(\xi)} =  \frac{e^{-(a+\kappa_+)\xi} u_*(\xi) + e^{-(a+\kappa_+)\xi}x(\xi)}{e^{-(a+\kappa_+)\xi} v_*(\xi) + e^{-(a+\kappa_+)\xi} y(\xi)}.$$ For $\xi\ra+\infty$, the nominator converges due to $u_*(\xi)\sim S_2 e^{\left(\frac{\chi}{\alpha}\kappa_+\right)\xi}$, (R1) and $w_+\geq -(a+\kappa_+)$. Again the denominator converges to a nonzero number by $y\in U_2$. Using (\ref{eq:upr}), the term $u'/v$ is treated in an analogous way with the help of Lemmas \ref{sec:key1} and \ref{sec:asym1}. As above, now the convergence of $a_4$, $a_5$, $a_6$ as $\xi\ra + \infty$ is a consequence of Lemma \ref{sec:fraccon}. 

Now we determine the explicit values of $a_i^+$ for $(x,y)=(0,0)$. Due to (\ref{eq:uv}) and Lemma \ref{sec:key1} we obtain $v_*'/v_*\ra a+\kappa_+$ and  $v_*'(1/v_*)'\ra -(a+\kappa_+)^2$. Using (\ref{eq:autoeo}) we further obtain that $$\frac{v_*''}{v_*} = \frac{-f(p_*) + (2a-\mu)q_* + a^2p_*}{p_*}\ra a^2 + (2a-\mu)\kappa_+ - \nu= (a+\kappa_+)^2$$ as $\xi\ra +\infty$. With the help of Lemma \ref{sec:asym1} and (R1) we deduce that  $u_*(\xi)/v_*(\xi)\sim C\, e^{(-(a+\kappa_+) + \frac{\chi}{\alpha}\kappa_+)\xi}\ra 0$ as $\xi\ra +\infty$. In the same way we obtain $u_*'/v_*\ra 0$ as $\xi\ra +\infty$, using (\ref{eq:Gasym}),  (\ref{eq:upr}) and Lemma \ref{sec:asym1}. Employing (\ref{eq:ev}) and Lemma \ref{sec:asym1}, the values for $a_i^+$ follow as stated.
\eprf
\begin{prop}\label{sec:Lbounded}
 \textsl{The operator $L_{(x,y)}:D\ra X$ is bounded for any $(x,y)\in \OO$.}
\end{prop}
\bprf The operators $L_1:D_1\ra X_1$ and $L_4:D_2\ra X_2$ are well-defined and bounded, since their continuous coefficients have limits at $\pm\infty$. The map $L_2:D_2\ra X_1$ is well-defined, since e.g. $\eta_+a_4h'' = (a_4/\eta_-)\eta h''$ has limits at $\pm\infty$, due to Lemma \ref{sec:asym2}. For $h\in D_2$ we have $$\|L_2h\|_{X_1} \leq \|a_4/\eta_-\|\|\eta h''\| +  \|a_5/\eta_-\|\|\eta h'\| + \|a_6/\eta_-\|\|\eta h\| \leq C \|h\|_{D_2},$$ and therefore $L_2$ is bounded. The map $L_3:X_1\ra X_2$ is well-defined (see the exponential rate of $a_7$ at $-\infty$ in Lemma \ref{sec:asym2}), and for $g\in X_1$ we have $\|L_3 g\|_{X_2} \leq \|\eta_-a_7\|\|\eta_+ g\| \leq C \|g\|_{X_1},$ hence $L_3$ is bounded.
\eprf
We now show that $F$ is Fr\'echet differentiable. For later purposes in Section \ref{sec:instability}, we show that the local approxmation by its linearization is better than one actually needs for differentiability. 

Let $E_1$, $E_2$ be two Banach spaces, $U\subset E_1$ be open and $f:U\ra E_2$. The map $f$ is called \textsl{$p$-linearizable} in $x_0\in U$ for some $p>1$, if there exists a bounded linear map $A:E_1\ra E_2$ such that 
\beq \label{eq:plin}
\|f(x_0+h)-f(x_0)-Ah\|_{E_2} = \OO(\|h\|_{E_1}^p) \qquad \trm{as }h\ra 0.
\eeq
Note that in this case $f$ is differentiable in $x_0$. If $f$ is $p$-linearizable in $x_0$ for $p>1$, then it is $q$-linearizable for any $q\in (1,p]$. We call $f$ $p$-linearizable, if it is $p$-linearizable in each $x_0\in U$. Taylor's formula implies that $C^2$ maps are $p$-linearizable for $p\in (1,2]$. 

Compositions of $p$-linearizable maps are again $p$-linearizable, which is proved similar to the chain rule. Thus we investigate the constituents of the first component of $F$. We will often use that for $f,g\in C^1(\ol{\R})$ we have $$\|fg\|\leq \|f\|\|g\|, \qquad \|fg\|_{C^1}\leq \|f\|_{C^1}\|g\|_{C^1}.$$ Define the weighted spaces 
\beq\label{eq:wsp1}
C_{\eta_+}(\ol{\R}) = \{x\in C(\ol{\R})\,|\, \eta_+x \in C(\ol{\R})\}, \quad \|x\|_{\eta_+}=\|\eta_+ x\|,
\eeq and 
\beq\label{eq:wsp2}
C_{\eta_+}^1(\ol{\R})=\{x\in C^1(\ol{\R})\,|\; \eta_+x, \eta_+x' \in C(\ol{\R})\},\quad \|x\|_{C_{\eta_+}^1}=\|\eta_+ x\| + \|\eta_+ x'\|.
\eeq The derivative $\partial_\xi:C_{\eta_+}^1(\ol{\R})\ra C_{\eta_+}(\ol{\R})$ is bounded and linear, and therefore $p$-linearizable for every $p>1$. The first summand in the first component of $F(x,y)$ equals
\beq\label{eq:fsfc}
d'(u_*+x)(u_*'+x')^2 + (\alpha+ d(u_*+x))(u_*'' +x'').
\eeq
\begin{lem}\label{sec:subst}
\textsl{For any function $f\in C^2(\R)$, the corresponding substitution operator $x\mapsto f\circ x$, $X_1\ra C(\ol{\R})$, is 2-linearizable with derivative $g\mapsto f'(x_0)\cdot g$ at $x_0\in X_1$.}
\end{lem}
\bprf Note that the substitution operator is well-defined and that the map $g\mapsto f'(x_0)\cdot g$ is bounded from $X_1$ to $C(\ol{\R})$. Let $x_0\in X_1$. Applying Taylor's formula pointwise, for any $\xi\in \R$ there is a $\tau=\tau(\xi)\in(0,1)$ such that $$\left | \left [ f(x_0+g) - f(x_0) - f'(x_0)\cdot g\right ](\xi)\right | = \left | \left [ \frac{1}{2} f''(x_0+\tau g) g^2 \right](\xi)\right|.$$ Taking the sup-norm, expanding the right-hand side by $\eta_+$ and using the local boundedness of $f''$ gives the $2$-linearizability.
\eprf
Note that to apply Lemma \ref{sec:subst} to (\ref{eq:fsfc}) we have to assume that $d\in C^3$.
\begin{lem}
\textsl{The multiplication, considered as a map $$C(\ol{\R})\times X_1\ra X_1 \quad \trm{or}\quad X_1\times X_1\ra X_1 \quad \trm{or} \quad X_1\times X_2\ra X_2,$$ is $2$-linearizable. In each case, the Fr\'echet derivative at $(x,y)$ is the map $(g,h)\mapsto xh + yg$, where $(x,y)$ and $(g,h)$ belong to the product spaces above.}
\end{lem}
\bprf
In each case, the multiplication is well-defined, and the stated linearization is bounded and linear. Take, for instance, $(x,y), (g,h)\in X_1\times X_2$. Then $$\|(x+g)(y+h) -xy  - xh-yg\|_{X_2} = \|gh\|_{X_2} \leq \|g\|_{X_1}\|h\|_{X_2}.$$ The other cases are treated analogously.
\eprf
We now treat the most difficult term of the first component of $F$.
\begin{lem}\label{sec:Q}
\textsl{The map $$Q: \OO\ra C_{\eta_+}^1(\ol{\R}), \qquad Q(x,y) = \frac{(u_*+x)(v_*+y)'}{v_*+y}$$ is $2$-linearizable with Fr\'echet derivative at $(x,y)\in \OO$ as stated in (\ref{eq:Qder}) below.}
\end{lem}
\bprf
We write $(u,v) = (u_*+x,v_*+y)$. The operator $Q$ maps to $C_{\eta_+}^1(\ol{\R})$, due to the Lemmas \ref{sec:asym1} and \ref{sec:fraccon}. The linearization of $Q$ in $(x,y)\in \OO$ is the map 
\beq\label{eq:Qder}
(g,h)\mapsto \frac{u}{v}\, h' + \frac{v'}{v}\, g - \frac{uv'}{v^2}\, h.
\eeq 
This is a bounded operator $D\ra C_{\eta_+}^1(\ol{\R})$, see Lemma \ref{sec:fraccon} and the proof of Lemma \ref{sec:asym2} for the fact that $u/(\eta v), u'/(\eta v)\in C(\ol{\R})$. Let $\|h\|_{D_2}$ be small enough such that $y+h\in U_2$. Then Lemmas \ref{sec:U2} and \ref{sec:fraccon} apply to $v+h$. We calculate
\begin{align*}
&\left \|(u+g)\frac{v'+h'}{v+h}- u\frac{v'}{v} - \frac{u}{v}\,h' - \frac{v'}{v}\,g + \frac{uv'}{v^2}\,h \right \|_{C_{\eta_+}^1} \\
& =  \left \| \frac{1}{v^2(v+h)}\left ( gh'v^2- uh'vh - gv'v h + uv' h^2\right) \right\|_{C_{\eta_+}^1} \\
& \leq  \left\|\frac{gh'}{v+h}\right\|_{C_{\eta_+}^1} +  \left\|\frac{uv'h^2}{v^2(v+h)}\right\|_{C_{\eta_+}^1} +  \left\|\frac{uh'h}{v(v+h)}\right\|_{C_{\eta_+}^1} +  \left\|\frac{gv'h}{v(v+h)}\right\|_{C_{\eta_+}^1}.
\end{align*}
The first summand is estimated by
$$\left\|\frac{gh'}{v+h}\right\|_{C_{\eta_+}^1}  \leq  \|g\|_{C_{\eta_+}^1} \left \| \frac{h'}{v+h}\right \|_{C^1}\leq C \|g\|_{D_1} \|h\|_{D_2}$$
since $$\left \|\frac{ h'}{v+h}\right \|\leq C \|h\|_{D_2}, \quad \left \| \left ( \frac{h'}{v+h}\right )' \right \| = \left \| \frac{h''}{v+h} - \frac{h'(v'+h')}{(v+h)^2}\right \|\leq C \|h\|_{D_2},$$
using Lemma \ref{sec:fraccon} and that, due to Lemma \ref{sec:U2}, the function $\eta (v+h)$ is uniformly bounded away from zero. 
The second term is estimated by
$$\left\|\frac{uv'h^2}{v^2(v+h)}\right\|_{C_{\eta_+}^1} \leq   \|\eta_+u\|_{C_{\eta_+}^1} \left \|\frac{v'}{v}\right\|_{C^1}\left\|\frac{h}{v} \right\|_{C^1}\left\|\frac{h}{v+h} \right\|_{C^1} \leq C \|h\|_{D_2}^2,$$  again due to Lemmas \ref{sec:U2} and \ref{sec:fraccon}. The third and forth summand are treated in a similiar fashion. We obtain that each summand is $\OO(\|(g,h)\|_D^2)$. 
\eprf
The map $F$ is composed of the derivative and maps treated in Lemmas \ref{sec:subst}-\ref{sec:Q}. Carefully composing the derivatives stated in these lemmas, we obtain the following result.
\begin{prop} \label{eq:Fplin}
 \textsl{The map $F:\OO\ra X$ is 2-linearizable with Fr\'echet derivative $F'(x,y) = L_{(x,y)}$ for $(x,y)\in \OO$.}
\end{prop}
By this proposition and the next one, $F$ fulfills (P1) and (P2).
\begin{prop}
 \textsl{The map $F':\OO\ra B(D,X)$ is locally Lipschitz continuous.}
\end{prop}
\bprf We have to show that $$\sup_{\|(g,h)\|_D=1} \|(F'(x_1,y_1) - F'(x_2,y_2))(g,h)\|_X \leq C\, \|(x_1-x_2,y_1-y_2)\|_D$$ locally holds in $\OO$. Recall that $F'(x,y) = L_{(x,y)}$ in (\ref{eq:A}). Expanding the coefficients of $L_2$ and $L_3$ by $\eta_-$, we obtain for $(x,y)\in \OO$ and $\|(g,h)\|_D = 1$
\begin{align}
&\|L_{(x,y)}(g,h)\|_{X} \nonumber \\
\leq &\max\{\trm{\,sup-norm of: } c, a_1, a_2, a_3, a_8, a_4/\eta_-, a_5/\eta_-, a_6/\eta_-, a_7/\eta_-\}.\label{eq:coeffLip}
\end{align}
See (\ref{eq:coeff}) for $a_1,...,a_8$. We assert that the coefficients in (\ref{eq:coeffLip}), considered as maps $\OO\ra C(\ol{\R})$, depend locally Lipschitzian on $(x,y)\in \OO$. All non-fraction terms depend locally Lipschitzian on $(x,y)$, using that the substitution operators and the multiplications are locally Lipschitzian. For the fractions we write $(u_1,v_1) = (u_* + x_1,v_*+y_1)$, $(u_2,v_2) = (u_* + x_2,v_*+y_2).$ Note that, for instance, $u_1-u_2=x_1-x_2$. For $v'/v$ we estimate
\begin{align*}
\left \|\frac{v_1'}{v_1} - \frac{v_2'}{v_2}\right \|& =  \left \| \frac{\eta(v_2(v_1'-v_2') + v_2'(v_2-v_1))}{\eta v_1v_2} \right \|\\
\leq  \frac{1}{\min|\eta v_1| } &\left \| \eta(y_1'-y_2')\right \| +\frac{1}{\min|\eta v_1 |} \cdot \left \| \frac{v_2'}{v_2}\right\|\cdot \left \| \eta(y_2-y_1)\right \|\leq C \|y_1-y_2\|_{D_2},
\end{align*}
using Lemmas \ref{sec:U2} and \ref{sec:fraccon}. Therefore also $v'(1/v')'= - (v'/v)^2$ is done. In a similiar way one treats $v''/v$. So the assertion holds for $a_2$ and $a_3$. For $u/(\eta_-v)$ we estimate
\begin{align*}
\left \|\frac{u_1}{\eta_-v_1} - \frac{u_2}{\eta_-v_2}\right \|   = &  \left \| \frac{\eta_+\eta(v_2(u_1-u_2) + u_2(v_2-v_1))}{\eta_-\eta_+\eta v_1v_2} \right \|\\ 
\leq &\; \frac{1}{\min|\eta v_1|} \left \| \eta_+(x_1-x_2)\right \| + \frac{\|\eta_+ u_2\|}{\min|\eta v_1\cdot\eta v_2 |} \left \| \eta(y_2-y_1)\right \|\\ 
\leq &\; C\left (\|x_1-x_2\|_{D_1} + \|y_1-y_2\|_{D_2} \right), 
\end{align*}
using Lemma \ref{sec:U2}. Similiarly one treats $u'/(\eta_-v)$. The assertion concerning $a_4/\eta_-$, $a_5/\eta_-$ and $a_6/\eta_-$ follows.
\eprf
To show sectoriality of $L_{(x,y)}$ for $(x,y)\in \OO$, we use the special structure of this operator, cf. (\ref{eq:A}). We exploit the fact that its  unbounded coefficients only occur in coupled terms of the first equation. Its main diagonal entries, $L_1$ and $L_4$, are sectorial by a standard result.
\begin{prop}
\textsl{For every $(x,y)\in \OO$ the operators $L_1= a_1\partial_{\xi\xi} + a_2 \partial_\xi + a_3$ and $L_4=  \partial_{\xi\xi} + c\partial_\xi + a_8$ on $X_1$ resp. $X_2$ with domains $D_1$ resp. $D_2$ are sectorial.}
\end{prop}
\bprf We first consider $L_4$. The maps $y\mapsto \eta y$, $X_2\ra C(\ol{\R})$, and $z\mapsto \eta^{-1}z$, $C^2(\ol{\R})\ra D_2$, are continuous isomorphisms. Thus $L_4$ is sectorial on $D_2$ if and only if $\wt{L_4}:=\eta L_4 \eta^{-1}$ is sectorial on $C(\ol{\R})$ with domain $C^2(\ol{\R})$, since the spectrum and resolvent estimates remain invariant under this similiarity transformation. For $h\in C^2(\ol{\R})$ one calculates $$\wt{L_4}h = h''+ (2\eta(\eta^{-1})' + c) h' + (a_8 + c\eta(\eta^{-1})' + \eta(\eta^{-1})'')  h.$$ Since $\eta$ is an exponential weight, the coefficents of $\wt{L_4}$ are continuous and have limits at $\pm \infty$. Now \cite[Theorem VI.4.3]{EN00} shows that $\wt{L_4}$ is sectorial. Noting that $a_1\in C^1(\ol{\R})$, $a_1\geq \alpha>0$, and $a_2,a_3\in C(\ol{\R})$ for any $(x,y)\in \OO$ due to Lemma \ref{sec:asym2}, one treats $L_1$ in a similiar fashion.
\eprf
The next lemma shows that the $\partial_{\xi\xi}$-bound of $\partial_\xi$ in the weighted space $X_2$ is zero.
\begin{lem}\label{sec:fde}
\textsl{For every $\eps>0$ there is $C_\eps>0$, such that for any $h\in D_2$ it holds that $$\|h'\|_{X_2}\leq \eps \|h''\|_{X_2}+C_\eps\|h\|_{X_2}.$$}
\end{lem}
\bprf Using that $\eta$ is an exponential weight, we obtain for arbitrary $\delta >0$ $$\|(\eta h)'\| \leq \delta\|(\eta h)''\| + C_\delta\|h\|_{X_2} \leq \delta\|h''\|_{X_2} + \delta C_0\|h'\|_{X_2} + C_\delta \|h\|_{X_2},$$ where we used (the proof of) \cite[Example III.2.2]{EN00} for the first inequality. The constant $C_\delta>0$ depends on $\delta$, but $C_0>0$ does not. Plugging this into $\|h'\|_{X_2}\leq \|(\eta h)'\| + C\|h\|_{X_2}$ and choosing $\delta$ small enough gives the result.
\eprf
\begin{lem}
\textsl{The operator $L_2= a_4\partial_{\xi\xi} + a_5 \partial_\xi + a_6: D_2\ra X_1$ is $L_4$-bounded.}
\end{lem}
\bprf As shown in the proof of Proposition \ref{sec:Lbounded}, $L_2$ is defined as stated. Using Lemma \ref{sec:fde}, we calculate for $h\in D_2$ 
\begin{align}
\|&L_2h\|_{X_1} \leq  \|a_4/\eta_-\| \|h''\|_{X_2} + \|a_5/\eta_-\| \|h'\|_{X_2} + \|a_6/\eta_-\| \|h\|_{X_2} \nonumber\\
&\leq C\|L_4h\|_{X_2} + (C\|c/\eta\|+\|a_5/\eta_-\|) \|h'\|_{X_2} + (C\|a_8/\eta\|+ \|a_6/\eta_-\|) \|h\|_{X_2}\nonumber \\
&\leq C\left (\|L_4h\|_{X_2} + \|h\|_{X_2}\right) + C\|h''\|_{X_2}.\label{eq:haha}
\end{align}
Using Lemma \ref{sec:fde} again, we estimate $\|h''\|_{X_2} \leq \|L_4 h\|_{X_2} +C \|h\|_{X_2} + \frac{1}{2} \|h''\|_{X_2}.$ Subtracting $\frac{1}{2} \|h''\|_{X_2}$ and plugging the result into (\ref{eq:haha}) finishes the proof.
\eprf
\begin{prop}
\textsl{For each $(x,y)\in \OO$, the operator $L_{(x,y)}$ is sectorial on $X$ with domain $D$.}
\end{prop}
\bprf Using that $L_1$, $L_4$ are sectorial, the $L_4$-boundedness of $L_2$, and that $L_3$ is a bounded operator, the result immediately follows from \cite[Corollary 3.3]{Nag89} and the bounded perturbation theorem for sectorial operators (\cite[Theorem III.2.10]{EN00}). 
\eprf
The sectoriality of $L_{(x,y)}$ for each $(x,y)\in \OO$ implies its closedness. Therefore its domain $D$, equipped with the graph norm $\|\cdot\|_{L_{(x,y)}}$ of $L_{(x,y)}$, is a Banach space. Estimating as in the proof of Proposition \ref{sec:Lbounded}, one verifies that the identity map $(D,\|\cdot\|_D)\ra (D, \|\cdot\|_{L_{(x,y)}})$ is bounded, and the open mapping theorem shows that the inverse identity map is also bounded. Thus $\|\cdot\|_{L_{(x,y)}}$ and $\|\cdot\|_D$ are equivalent norms on $D$ for each $(x,y)\in \OO$. 

This finally shows that $F$ enjoys the properties (P1)-(P3) from Theorem \ref{sec:LW}, and Theorem \ref{sec:thmlw} is proved.

\section{Proof of Theorem \ref{sec:spec}: A Part of the Spectrum of $F'(0,0)$}\label{spectrum}
In this section we show that under certain restrictions on the coefficients of the model, the wave speed and the weight on the right half-line (see (R2) in Section \ref{sectionfr}), the operator $A=F'(0,0)=L_{(0,0)}$ has spectral values with positive real part. Recall that $A= M_1\partial_{\xi\xi}+ M_2\partial_\xi + M_3,$ where $$M_1= \left (\begin{array}{cc}a_1 & a_4 \\ 0 & 1 \end{array} \right ), \quad M_2= \left (\begin{array}{cc}a_2 & a_5 \\ 0 & c \end{array} \right ), \quad M_3 = \left (\begin{array}{cc}a_3 & a_6 \\ a_7 & a_8 \end{array} \right ),$$ is considered as an operator on $X$ with domain $D$. See (\ref{eq:coeff}) with $(u,v)$ replaced by $(u_*,v_*)$ for the coefficents $a_i$, and (\ref{eq:apmlimits}) for their limits at $+\infty$.
 
Being aware of the different notions for ``essential spectrum'' in the literature (cf. \cite[Section IV.5.6]{Kat66}), we define $$\sigma_{ess}(A)=\{\lambda\in \C\;|\; A_\lambda=A-\lambda \trm{ is not a Fredholm operator}\}.$$ The essential fact for $\sigma_{ess}$ is that for closed operators it remains invariant under relatively compact perturbations (\cite[Theorem IV.5.26]{Kat66}). In Section 3 of the survey \cite{San02}, a machinery for calculating $\sigma_{ess}$ for second order ordinary differential operators with coefficents in $C(\ol{\R})$ is described. Since the coefficients of $A$ are not bounded on $\R_-$, see Lemma \ref{sec:asym2}, we will remove $\R_-$ and apply this machinery for coefficents in $C(\ol{\R}_+)$, where this space and $C_{\eta_+}^2(\ol{\R}_+)$ are defined analogously to (\ref{eq:wsp1}) and (\ref{eq:wsp2}).
\begin{prop}\label{sec:trunc}
 \textsl{Suppose $A_\lambda= M_1\partial_{\xi\xi}+ M_2\partial_\xi + (M_3-\lambda)$ is a Fredholm operator on $X$ with domain $D$. Then $B_\lambda$, defined as $A_\lambda$ but on $X_+ = C_{\eta_+}(\ol{\R}_+)^2$ with domain $D_+ = C_{\eta_+}^2(\ol{\R}_+)^2$, is a Fredholm operator as well.}
\end{prop}
\bprf
We show that $\dim (\ker B_\lambda)<+\infty$ and $\trm{codim}(\im B_\lambda) <+\infty$. The map $f\mapsto M_1^{-1}f$, $X_+\ra X_+$, is a continuous isomorphism. Now $M_1^{-1}B_\lambda$ is a second order ordinary differential operator on $\R_+$, thus its kernel is finite dimensional. We conclude that $\ker B_\lambda$ is finite dimensional. Since $A_\lambda$ is supposed to be Fredholm, we have $X = (\im A_\lambda)\oplus V$ with $\dim V<+\infty$. Let $f\in X_+$. Extend $f$ to a function $\wt{f}\in X$, then $\wt{f} = A_\lambda \wt{u} + \wt{v}$ for some $\wt{u}\in D$ and $\wt{v} \in V$. Restricting $\wt{u}$ and $\wt{v}$ on $\R_+$ to functions $u\in D_+$ and $v\in V|_{\R_+}$, we obtain $f= B_\lambda u + v$ on $\R_+$, and therefore $X_+ =  (\im B_\lambda) + V|_{\R_+}$. Since $V|_{\R_+}$ is finite dimensional, we obtain $X_+ = (\im B_\lambda) \oplus W$, where $W$ is a complement of $(\im B_\lambda\cap V|_{\R_+})$ in $V|_{\R_+}$. Thus $\im B_\lambda$ has a finite dimensional co-image. 
% $D_+ = (\im A_\lambda)|_{\R_+}\oplus W$ with $W \cong V|_{\R_+}/((\im A_\lambda)|_{\R_+}\cap V|_{\R_+})$. By restricting again it holds that $(\im A_\lambda)|_{\R_+} \subset \im B_\lambda$. Setting $U= \im B_\lambda \cap W$, we obtain $D_+ = (\im A_\lambda)|_{\R_+}\oplus U \oplus Z$ with $Z\cong W/U$, and $\im B_\lambda = (\im A_\lambda)|_{\R_+}\oplus U$. This shows that $B_\lambda$ has a finite dimensional co-image. 
\eprf
Now consider the continuous isomorphisms $$f\mapsto \eta_+ f,\quad  X_+\ra C(\ol{\R}_+); \qquad g \mapsto \eta_+^{-1}g, \quad  C^2(\ol{\R}_+)\ra D_+.$$ The operator $B_\lambda$ is Fredholm if and only if $\wt{B_\lambda}= \eta_+ B_\lambda \eta_+^{-1}:C^2(\ol{\R}_+)^2\ra C(\ol{\R}_+)^2$ is Fredholm. We calculate $\wt{B_\lambda} = \wt{M_1}\partial_{\xi\xi} + \wt{M_2}\partial_\xi + \wt{M_3}$ with the coefficients $$\wt{M_1} = M_1 = \left ( \begin{array}{cc} a_1 & a_4 \\ 0 & 1 \end{array} \right ),\qquad \wt{M_2} = \left ( \begin{array}{cc} -2 w_+ a_1 +a_2 & -2w_+a_4 + a_5 \\ 0 & -2w_+ + c \end{array} \right ),$$ 
\beq\label{eq:MSchlange}
\wt{M_3} = \left ( \begin{array}{cc} w_+^2a_1 - w_+a_2 + a_3 - \lambda & w_+^2a_4 - w_+a_5 + a_6 \\ a_7 & w_+^2- w_+ c + a_8 - \lambda \end{array} \right ).
\eeq The multiplication by $M_1$ is an isomorphism as a map $C(\ol{\R}_+)^2\ra C(\ol{\R}_+)^2$ (note that $a_1\geq \alpha>0$ and that all coefficents are bounded on $\R_+$), hence also $M_1^{-1}\wt{B_\lambda}:C^2(\ol{\R}_+)^2\ra C(\ol{\R}_+)^2$ is a Fredholm operator. Set $C_0(\R_+)=\{f\in C(\R_+)\,|\,f(+\infty) = 0\}$ and $C_0^2(\R_+)= \{f\in C^2(\R_+)\,|\,f,f',f''\in C_0(\R_+)\}.$ Then $$C(\ol{\R}_+) = C_0(\R_+)\oplus \R, \qquad C^2(\ol{\R}_+)= C_0^2(\R_+)\oplus \R,$$ since $f = (f-f(+\infty)) + f(+\infty)$ and necessarily $f'(+\infty) = f''(+\infty) = 0$ for $f\in C^2(\ol{\R}_+)$. The operator $M_1^{-1}\wt{B_\lambda}$ maps elements of $C_0^2(\R_+)^2$ into $C_0(\R_+)^2.$  Thus we can define $$D_\lambda= \partial_{\xi\xi} + M_1^{-1}\wt{M_2}\partial_\xi + M_1^{-1}\wt{M_3}: C_0^2(\R_+)^2\ra C_0(\R_+)^2,$$ and $D_\lambda$ is Fredholm if and only if $M^{-1}\wt{B_\lambda}$ is Fredholm. The corresponding first order operator of $D_\lambda$ is $E_\lambda= \partial_\xi+T_\lambda:C_0^1(\R_+)^4\ra C_0(\R_+)^4$, where 
\beq \label{eq:tlambda}
T_\lambda = \left ( \begin{array}{cc} 0 & -\id \\ M_1^{-1}\wt{M_3} & M_1^{-1}\wt{M_2}\end{array}\right ).
\eeq Setting $T_\lambda^+ = T_\lambda(+\infty)$, it is a well-known fact that $D_\lambda$ is a Fredholm operator if and only if $T_\lambda^+$ is a hyperbolic matrix (see \cite[Chapter 3]{San02} for the $L^2$ and $C_b$ case, and \cite[Appendix to Chapter 5]{Hen81}). 

For completeness, we sketch the proof for the $C_0(\R_+)$ case, following \cite{San02}. The Fredholm properties of $D_\lambda$ are the same as the Fredholm properties of $E_\lambda$ (\cite[Appendix A]{SS08}; the proof there is easily adopted to the $C_0(\R_+)$ case). Replacing $T_\lambda$ by $T_\lambda^+$, we receive the constant coefficient operator $E_\lambda^+ = \partial_\xi + T_\lambda^+$, which differs from $E_\lambda$ by a relatively compact perturbation (\cite[Appendix to Chapter 5]{Hen81}). Therefore $E_\lambda$ and $E_\lambda^+$ have the same Fredholm properties (\cite[Theorem IV.5.26]{Kat66}). It follows from the proof of \cite[Theorem 1]{Pal88}, that if $E_\lambda^+$ is a Fredholm operator, then the corresponding homogenous equation $u' + T_\lambda^+ u=0$ possesses an exponential dichotomy (for a definition and properties see \cite{Cop78}; the converse of this statement is also true, see the proof of \cite[Lemma 4.2]{Pal84}). Being precise, \cite[ Lemma 1]{Pal88} shows that $E_\lambda^+$ always has dense image. Thus if $E_\lambda^+$ is supposed to be Fredholm, it must be surjective. As indicated in \cite{Pal88}, now it follows from \cite[Theorem 64B]{MS66} that the corresponding homogenous equation possesses an exponential dichotomy. From \cite[Chapter 6]{Cop78} it follows that $u'+T_\lambda^+u=0$ possesses an exponential dichotomy if and only if $T_\lambda^+$ is a hyperbolic matrix.

\begin{rem} We emphasize that this also shows that the operator $B_\lambda$ is surjective if and only if $T_\lambda^+$ is hyperbolic.
\end{rem}

Summarizing, if $T_\lambda^+$ possesses a purely imaginary eigenvalue then $A_\lambda$ is not a Fredholm operator, and therefore $\lambda\in \sigma_{ess}(A)$. Thus for $\lambda\in \C$ we are looking for solutions $h\in \R$ of the so-called dispersion relation $\det(T_\lambda^+ -\ii h) = 0$. Substituting the limits (\ref{eq:apmlimits}) into (\ref{eq:MSchlange}) we obtain from (\ref{eq:tlambda}) $$T_\lambda^+ = \left ( \begin{array}{cccc} 0 & 0 & -1 & 0 \\ 0 & 0 & 0 & -1 \\ w_+^2 - \frac{w_+a_2^+ + \lambda}{\alpha} & 0 & -2w_+ + \frac{a_2^+}{\alpha} & 0  \\ 0 & w_+^2 - w_+c + \gamma- \lambda & 0 & -2w_++c  \end{array} \right ),$$
where $a_2^+ =  c + \frac{\chi}{2}\left (c - \sqrt{c^2-4 \gamma}\right)$. We calculate 
\bea
\det(T_\lambda^+ -\ii h)& = & (-h^2 + \ii(2w_+- a_2^+/\alpha)h + w_+^2 - w_+a_2^+/\alpha - \lambda/\alpha) \\
& &  \cdot\, (-h^2 + \ii( 2w_+-c) h + w_+^2 - w_+c + \gamma -\lambda),
\eea
and obtain the following result.
\begin{prop}
\textsl{The sets 
\begin{align*}
S_1 = & \left \{\lambda \in \C\,|\,\emph{\rea} \lambda = \alpha(-h^2 + w_+^2 - w_+ a_2^+/\alpha), \;\emph{\ima} \lambda = (2\alpha w_+ - a_2^+)h, h\in \R\right\},\\
S_2 = & \left \{\lambda \in \C\,|\,\emph{\rea} \lambda = -h^2 + w_+^2 - w_+c +\gamma,\;\emph{\ima} \lambda = (2w_+-c)h, h\in \R\right \},
\end{align*}
are contained in $\sigma(A)$.}
\end{prop}

The sets $S_1$, $S_2$ are shaped as parabolas, open to the left. We are interested in spectral values with positive real parts. For the weight we allowed values $w_+\in [-(a+\kappa_+), -\frac{\chi}{\alpha}\kappa_+]$. One easily checks that $w_+^2 - w_+ a_2^+/\alpha\leq 0$ for every choice of $w_+$, so in $S_1$ numbers with positive real part do not occur. 

For $S_2$, consider the polynomial $t^2-ct + \gamma$. Its roots are $t_{\pm} = c/2 \pm \sqrt{c^2-4\gamma}/2>0$, i.e. for $t<t_-$ and $t>t_+$ the polynomial has positive values. We have $-(a+\kappa_+)= t_-$ (see (\ref{eq:ev})). If the  inequality $-\frac{\chi}{\alpha}\kappa_+> t_+$ holds, then we can always choose an exponential rate $w_+\in \left (t_+,-\frac{\chi}{\alpha}\kappa_+\right ]$ on $\R_+$, such that spectral values with positive real part occur. The last inequality is equivalent to $c(\chi-\alpha+2)> (\alpha+\chi) \sqrt{c^2-4\gamma}.$ This can never hold if $\chi\leq \alpha -2$. In the case $\chi >\alpha-2$, we can rewrite it to $c^2(\chi+1)(\alpha-1)< \gamma(\chi+\alpha)^2$, which always holds for $\alpha\leq 1$. For $\alpha>1$ this yields the upper bound (\ref{eq:crest}) on the wave speed. 

These restrictions are summarized in (R2) in Section \ref{sectionfr}, and thus Theorem \ref{sec:spec} is proved.

\begin{rem}\label{remtrunk} The proof of Proposition \ref{sec:trunc} shows that by truncating $\R$ to a half-line, the Fredholm index of $A_\lambda$ will increase in general. This is the reason why we did not choose the essential spectrum to be the (larger) set of $\lambda\in \C$, for which $A_\lambda$ is not a Fredholm operator of index zero, as it is done in \cite{San02}. 

Further, the converse of this proposition is wrong. Assume the converse was true. Choose matrices $T_-$, $T_+$ of type (\ref{eq:tlambda}), where $T_+$ is hyperbolic but $T_-$ is not, and a matrix-valued function $T\in C(\ol{\R})$ with $T(\pm\infty) = T_\pm$. Now consider the second order operator $S$ which corresponds to $\partial_\xi + T$ as $S:C_0^2(\ol{\R}_+)\ra C_0(\ol{\R}_+)$. By assumption, $S$, considered on $\R$, is Fredholm. But applying Proposition \ref{sec:trunc} and the machinery described above on $\R_-$ gives a contradiction. 

This means that by truncating $\R_-$ in Proposition \ref{sec:trunc}, we are in general only able to find a part of the essential spectrum of the original operator. Therefore, determining stability of a (family of) steady state(s) is in general not possible when using Proposition \ref{sec:trunc}.

Since the travelling wave problem is invariant under translations, we obtain a trivial zero eigenvalue of $F'(0,0)$ if $(u_*,v_*)\in \OO$, i.e. if the wave itself is contained in the space of perturbations (see Remark \ref{sec:reminst} and Lemma \ref{sec:asym1}). However, this fact is of no interest to us in the present situation.
\end{rem}

\section{Proof of Theorem \ref{sec:linin}: Instability Without Spectral Gap}\label{sec:instability}
In an abstract setting, we show the principle of linearized instability for fully nonlinear parabolic problems, without assuming the existence of a spectral gap of the linearization in a steady state. This is the case for a travelling wave from Theorem \ref{sec:ex}, see Theorem \ref{sec:spec}.

Suppose $X,D,\OO$ and $F:\OO\ra X$ are as in Theorem \ref{sec:LW}, such that $F$ enjoys properties (P1)-(P3). Then the evolution equation $\uu_t=F(\uu)$ is locally well-posed in the sense (LW1), (LW2). Suppose further that $F$ is $p$-linearizable in $\uu_0=0$ for some $p>1$ (see (\ref{eq:plin})), which will be fixed from now on, and that the sectorial operator $$A=F'(0)$$ has a spectral value with positive real part. To prove nonlinear instability of $\uu_0=0$, we use the following result.

\begin{thm}[{\cite[Theorem 5.1.5]{Hen81}}]\label{sec:henry} \textsl{Suppose $D$ is a real Banach space and $U\subset D$ is an open neighbourhood of the origin. The map $T:U\ra D$ is supposed to be continuous with $T(0)=0$, and to be $q$-linearizable in zero for some $q>1$ by $M\in B(D)$ with spectral radius greater than one.}

\textsl{Then $\ue=0$ is unstable in the sense of Lyapunov, i.e. there is $\eps_0>0$ and in any neighbourhood of the origin in $U$ there is a $\ue_0$, such that for some natural number $N  \geq 1$, the sequence $(\ue_n)_{n=0,...,N}$, given by $\ue_n=T(\ue_{n-1})$, is defined and $\|\ue_N\|_D \geq \eps_0$.}
\end{thm}
Choose an arbitrary $\theta\in (0,1)$. Thanks to (LW2) there is an open set $U\subset \OO$ containing zero, such that $\tau(\uu_0)\geq 1$ for any $\uu_0\in U$, i.e. the solution map $\uu(\cdot ,\uu_0)$ is defined for $t\in [0,1]$. It is further locally Lipschitz continuous as a map from $U$ into $C_\theta^\theta(]0,1],D)$. Now the time-one map $$T=\uu(1,\cdot):U\ra D$$ is defined. To prove Theorem \ref{sec:linin} we show that $T$ satisfies the assumptions of Theorem \ref{sec:henry}, with Fr\'echet derivative $e^A:D\ra D$.

Since $F(0) = 0$, we have $T(0)=0$. Further, $T$ is locally Lipschitz continuous. Using the graph norm of $A$ on $D$ and \cite[Proposition 2.1.1]{Lun95}, one checks that $e^A:D\ra D$ is continuous. Since the spectral mapping theorem $\sigma(e^A) \bs\{0\}= e^{\sigma(A)}$ holds for a sectorial operator $A$ (\cite[Corollary 2.3.7]{Lun95}), $e^A$ has spectral radius greater than one, considered as a continuous operator on $X$. Since $A$ and $e^A$ commute on $D$ (\cite[Proposition 2.1.1]{Lun95}), we have $e^A\uu_0 = (\lambda-A)^{-1} e^A(\lambda - A)\uu_0$ for each $\uu_0\in D$ and arbitrary $\lambda$ in the resolvent set of $A$. Thus also $e^A$, considered as an element of $B(D)$, has spectral radius greater than one.

It remains to show the $q$-linearizability of $T$ in zero for some $q>1$. 

\begin{prop}\label{sec:keyprop}
\textsl{Suppose $F$ is $p$-linearizable in $\ue_0=0$ for some $p>1$. Then the time-one map $T:U\ra D$ is $q$-linearizable in $\ue_0=0$ for any $q\in (1,p)$ with Fr\'echet derivative $e^A\in B(D)$.}
\end{prop}
\bprf
Set $G=F-A$, then $G:\OO\ra X$ is continuous. For any $\uu_0\in U$ the corresponding solution $\uu(\cdot,\uu_0)$ fulfills $\uu_t = A\uu + G(\uu)$ for $t\in [0,1]$. Since $G(\uu(\cdot,\uu_0))\in C([0,1];X)$, the time-one map can be represented by the variation of constants formula (\cite[Proposition  4.1.2]{Lun95}): $$T(\uu_0) = e^A \uu_0 + \int_0^1 e^{(1-s)A} G(\uu(s,\uu_0))\D s$$ Consider the integral term as a map $R(\uu_0)$ for $\uu_0\in U$. Fix $q\in(1,p)$. We have to show that there is a $\delta >0$ such that $\|R(\uu_0)\|_D\leq C\,\|\uu_0\|_D^q$ for each $\|\uu_0\|_D\leq\delta$.

The family $\left ( e^{(1-s)A}\right)_{s\in[0,1]} \subset B(X)$ is uniformly bounded (\cite[Proposition 2.1.1]{Lun95}). The local Lipschitz continuity of the solution map in $\uu_0=0$ and $F(0)=0$ imply that there is $\delta>0$ such that
\beq\label{eq:loclip}
\|\uu(\cdot,\uu_0)\|_{C_\theta^\theta(]0,1]; D)}\leq C\;\|\uu_0\|_D
\eeq if $\|\uu_0\|_D\leq \delta$. This yields in particular $\sup_{s\in [0,1]} \|\uu(s,\uu_0)\|_D\leq C\delta$. Choosing $\delta$ small enough, the $p$-linearizability of $F$ in zero  and (\ref{eq:loclip}) yield for any $s\in[0,1]$
\beq\label{eq:G}
\|G(\uu(s,\uu_0))\|_X\leq C\,\|\uu(s,\uu_0)\|_D^p \leq C\, \|\uu_0\|_D^p,
\eeq 
and therefore
%\label{eq:Rest}
$$\|R(\uu_0)\|_X \leq C\, \sup_{s\in[0,1]} \|G(\uu(s,\uu_0))\|_X\leq  C\, \|\uu_0\|_D^q,$$
provided $\|\uu_0\|_D\leq\delta$.

By making $\delta$ once more smaller if necessary, due to Lemma \ref{sec:alpha} we have $G(\uu(\cdot,\uu_0))\in C_{\beta}^{\beta}(]0,1];X)$ and estimate (\ref{eq:n7}) for some $\beta\in (0,1)$ with $\beta<\theta$. \cite[Theorem 4.3.5]{Lun95} gives $R(\uu_0)\in D$ and $$\|AR(\uu_0)\|_X\leq C\; \|G(\uu(\cdot,\uu_0))\|_{C_{\beta}^{\beta}(]0,1];X)},$$ therefore (\ref{eq:n7}) finishes the proof.
\eprf
\begin{lem}\label{sec:alpha}
\textsl{In the setting of the proof above, there are $\delta >0$ and $\beta\in (0,1)$ with $\beta< \theta$ such that for $\|\uu_0\|_D\leq\delta$ we have $G(\uu(\cdot,\uu_0))\in C_\beta^\beta(]0,1];X)$ and the estimate
\beq\label{eq:n7}
\|G(\uu(\cdot,\uu_0))\|_{C_\beta^\beta(]0,1];X)} \leq C\, \|\uu_0\|_D^q.
\eeq}
\end{lem}
\bprf Fix $\uu_0\in U$ with $\|\uu_0\|_D\leq\delta$, where $\delta$ will be chosen small enough in the sequel. For simplicity we write $\uu(\cdot) = \uu(\cdot,\uu_0)$. Take $\beta\in (0,1)$ smaller than $\theta$ such that 
\beq\label{eq:choice}
p\left (1-\frac{\beta}{\theta}\right) + \frac{\beta}{\theta}>q.
\eeq
We have to show that $$\sup_{t\in[0,1]}\|G(\uu(t))\|_X + \sup_{\eps\in(0,1)} \eps^\beta [G(\uu(\cdot))]_{C^\beta([\eps,1];X)} \leq C \,\|\uu_0\|_D^q.$$ For the first summand see (\ref{eq:G}). In (P2), the derivative $F'$ is assumed to be locally bounded, thus $F:\OO\ra X$ is locally Lipschitz continuous near zero.  Therefore $G=F-A:\OO\ra X$ is Lipschitz continuous on a ball $B_{\delta_0}(0)$ for some  $\delta_0>0$.

If $\|\uu_0\|_D\leq\delta$ with $\delta$ small enough then $\sup_{t\in[0,1]} \|\uu(\cdot)\|_D\leq\delta_0$, thanks to (\ref{eq:loclip}). In this case we can estimate
\begin{align*}
& [G(\uu(\cdot))]_{C^\theta([\eps,1];X)}  =  \sup_{t,s\in[\eps,1],t>s} \frac{\|G(\uu(t))- G(\uu(s))\|_X}{(t-s)^\theta} \\
& \leq  C\sup_{t,s\in[\eps,1],t>s} \frac{\|\uu(t)- \uu(s)\|_D}{(t-s)^\theta} =  C\;[\uu(\cdot)]_{C^\theta([\eps,1];D)}
\end{align*}
for any $\eps\in (0,1)$, where the constant $C$ is independent of $\eps$. Hence, using (\ref{eq:loclip}) again,
\beq\label{eq:GHol}
\sup_{\eps\in(0,1)} \eps^\theta [G(\uu(\cdot))]_{C^\theta([\eps,1];X)} \leq C\; \|\uu(\cdot)\|_{C_\theta^\theta(]0,1];D)} \leq C\; \|\uu_0\|_D.
\eeq 
Now we claim that for $0<\beta<\theta <1$ there is a constant $K>0$, such that for any closed interval $I\subset \R$ it holds that
\beq\label{eq:Hint}
\|f\|_{C^\beta(I,X)}\leq K\; \|f\|_\infty^{1-\frac{\beta}{\theta}} \|f\|_{C^\theta(I,X)}^{\frac{\beta}{\theta}}\qquad \trm{ for any } f\in C^\theta(I,X).
\eeq In \cite[Proposition 1.1.3]{Lun95} the case $I=\R$ is treated. The general case follows from the fact that the H\"older norm of a function $f:I\ra \R$ is not changed if $f$ is constantly continued from $I$ to $\R$. Using (\ref{eq:Hint}), (\ref{eq:G}) and (\ref{eq:GHol}), we perform the following estimates:
\begin{align*}
& \sup_{\eps\in (0,1)}  \eps^\beta [G(\uu(\cdot))]_{C^\beta([\eps,1];X)}  \leq \sup_{\eps\in (0,1)} \eps^\beta \|G(\uu(\cdot))\|_{C^\beta([\eps,1];X)} \\
& \leq  K\; \|G(\uu(\cdot))\|_\infty^{1-\frac{\beta}{\theta}}\sup_{\eps\in (0,1)} \eps^\beta  \|G(\uu(\cdot))\|_{C^\theta([\eps,1];X)}^\frac{\beta}{\theta} \\
& \leq  K \; \|G(\uu(\cdot))\|_\infty^{1-\frac{\beta}{\theta}}  \left ( \|G(\uu(\cdot))\|_{\infty} + \sup_{\eps\in (0,1)} \eps^\theta [G(\uu(\cdot))]_{C^\theta([\eps,1];X)}\right )^{\frac{\beta}{\theta}}\\
& \leq  C\,\|\uu_0\|_D^{p\left (1-\frac{\beta}{\theta}\right )} \left ( C\|\uu_0\|_D^p +  C\|\uu_0\|_D\right )^{\frac{\beta}{\theta}} \leq C\;\|\uu_0\|_D^{p\left (1-\frac{\beta}{\theta}\right) + \frac{\beta}{\theta}}
\end{align*}
Our choice of the exponent in (\ref{eq:choice}) gives the result for $q\in (1,p)$.
\eprf
Thanks to Proposition \ref{sec:keyprop}, the time-one map $T$ generated by $F$ fulfills the assumptions of Theorem \ref{sec:henry}, and this proves Theorem \ref{sec:linin}.

\end{document}